\documentclass[11pt, reqno]{amsart}

 \usepackage{amsmath}
 \usepackage{amssymb}
\usepackage{bm}

\DeclareMathAccent{\mathring}{\mathalpha}{operators}{"17}

\newcommand{\mysection}[1]{\section{#1}
      \setcounter{equation}{0}}

\newcommand{\nlimsup}{\operatornamewithlimits{\overline{lim}}}

\newtheorem{theorem}{Theorem}[section]
\newtheorem{lemma}[theorem]{Lemma}

\newtheorem{corollary}[theorem]{Corollary} 

\theoremstyle{definition}
\newtheorem{assumption}{Assumption}[section]

\theoremstyle{remark}
\newtheorem{remark}{Remark}[section]

\newtheorem{example}{Example}[section]

\newcommand{\tr}{\text{\rm tr}\,}
\newcommand{\loc}{\text{\rm loc}}

 \makeatletter 
 \def\dashint{%
 \operatorname%
 {\,\,\text{\bf--}\kern-.98em\DOTSI\intop\ilimits@\!\!}}
\def\ninf{\qopname\relax\@empty{inf\phantom{p}\!\!\!}}
 \makeatother

 \newcommand{\WO}{\overset{ 
 \scriptscriptstyle0}%
 { W}\,\!}
 
\newcommand\gb{\mathfrak{b}}

\newcommand\bbeta{\text{\raise-.2ex\hbox{$\bm{\beta}$}}}

\newcommand\esssup{\operatornamewithlimits{esssup}}

\newcommand\bR{\mathbb{R}}

\newcommand\bS{\mathbb{S}}

\newcommand\cF{\mathcal{F}}

\newcommand\cL{\mathcal{L}}
\newcommand\cM{\mathcal{M}}

\newcommand\cN{\mathcal{N}}

\newcommand\frA{\mathfrak{A}}

\newcommand\frM{\mathfrak{M}}

\newcommand\dist{{\rm dist}\,}

\begin{document}

\title[Stochastic equations with drift in $L_{d}$]
{On stochastic equations with drift in $L_{d}$}

\author{N.V. Krylov}
 
\email{nkrylov@umn.edu}
\address{127 Vincent Hall, University of Minnesota,
 Minneapolis, MN, 55455}



\begin{abstract}
   For the It\^o stochastic equations
in $\bR^{d}$ with drift in $L_{d}$ several results
are discussed such as the existence of weak solutions,
the existence of the corresponding Markov process, the
Aleksandrov type estimates of their Green's functions,
which yield their summability to the power of $d/(d-1)$, the
Fabes-Stroock type estimates which show that 
Green's functions are summable to a higher degree, the
Fanghua Lin type estimates, which are one of the
 main tools in the $W^{2}_{p}$-theory of fully nonlinear
elliptic equations, the fact that Green's functions
are in the class $A_{\infty}$ of Muckenhoupt and a 
few other results.
\end{abstract}

\maketitle

\mysection{Introduction}

Let $\bR^{d}$ be an Euclidean space of points
$x=(x^{1},...,x^{d})$. We assume that $d\geq2$ and denote
$$
B_{R}(x)=\{y\in\bR^{d}:|y-x|<R\},\quad
B_{R} =B_{R}(0),\quad D_{i}=\frac{\partial}{\partial x^{i}},
\quad D_{ij}=D_{i}D_{j}.
$$

First we address the issue of   existence
of solutions of It\^o's equations with drift term in $L_{d}$.

\begin{example}
                                             \label{example 9.6.1}
Let  $b(x)=- x|x|^{-2}(d/2)$. We have $b\in L_{d-\varepsilon}(B_{1})$
for any $\varepsilon\in(0,1)$ but not for $\varepsilon=0$ and it turns
out that there is no solutions of the equation $dx_{t}=dw_{t}
+b(x_{t})\,dt$ starting at zero, where $w_{t}$ is a $d$-dimensional
Wiener process.

Indeed, if we assume the contrary, then for the equation to make sense
$$
\int_{0}^{t}|b(x_{t})|\,dt
$$
should be finite. On the other hand, by It\^o's formula
$|x_{t}|^{2}$ turns out to be at least a
 local martingale and, since it is nonnegative
and starts from zero, it is zero. Then
$$
0=x_{t}=w_{t}+\int_{0}^{t}b(x_{t})\,dt,
$$
that is the sum of 
$w_{t}$ and a function of bounded variation is zero,
which is impossible.
\end{example}

Thus, in the general case we can only hope that the existence
of solutions of stochastic equations  holds if $b\in L_{d}$.
To formulate a result which  contains the 
existence theorem we introduce some necessary objects.

Introduce $\bS$ as the set of $d\times d$ symmetric matrices,
and for $\delta\in(0,1)$ let $\bS_{\delta}$ be the subset
of $\bS$ consisting of matrices whose eigenvalues are between
$\delta$ and $\delta^{-1}$.

Let $b (x)$, $b^{(k)}(x)$, $k= 1,2,...$, be   $\bR^{d}$-valued 
Borel functions
on $\bR^{d}$ such that, for a constant $\|b\|<\infty$,
$$
\|b \|_{L_{d}(\bR^{d})},\|b^{(k)}\|_{L_{d}(\bR^{d})}\leq \|b\|,
\quad k= 1,2,...,
$$
and $b^{(k)}\to b $ as $k\to\infty$ in $L_{d}(\bR^{d})$.
 Let $a (x)$, $a^{(k)}(x)$, $k= 1,2,...$,    
be  Borel functions on $\bR^{d}$ with values in $\bS_{\delta}$
for some $\delta\in(0,1)$ such that
$a^{(k)}\to a $ as $k\to\infty$ (a.e.).

\begin{theorem}
                                           \label{theorem 9.6.4}
Take $x\in\bR^{d}$. 
(i) There exists 
a probability space $(\Omega ,\cF ,P )$,
a filtration of $\sigma$-fields $\cF _{t}\subset \cF $, $t\geq0$,
a process $w _{t}$, $t\geq0$, which is a $d$-dimensional Wiener process
relative to $\{\cF _{t}\}$, and an $\cF _{t}$-adapted
process $x_{t}$ such that 
 (a.s.) for all   $t\geq0$
\begin{equation}
                                                 \label{11.29.2}
x _{t}=x  +\int_{0}^{t}\sqrt{a (x_{s})}\,dw_{s}
+\int_{0}^{t}b (x_{s}) \,ds.
\end{equation}

(ii) Furthermore,
let $x^{(k)}\in\bR^{d}$, $k= 1,2,...$, and let $x^{(k)}
\to x $  as $k\to\infty$. 
Assume that for each $k=1,2,...$  
 there exists 
a probability space $(\Omega^{(k)},\cF^{(k)},\\P^{(k)})$,
a filtration of $\sigma$-fields $\cF^{(k)}_{t}\subset \cF^{(k)}$, $t\geq0$,
a process $w^{(k)}_{t}$, $t\geq0$, which is a $d$-dimensional Wiener process
relative to $\{\cF^{(k)}_{t}\}$, and an $\cF^{(k)}_{t}$-adapted
process $x^{(k)}_{t}$ such that (a.s.) for all   $t\geq0$
\begin{equation}
                                                 \label{11.29.1}
x^{(k)}_{t}=x^{(k)} +\int_{0}^{t}\sqrt{a^{(k)}(x^{(k)}_{s})}\,dw^{(k)}_{s}
+\int_{0}^{t}b^{(k)}(x^{(k)}_{s}) \,ds.
\end{equation}

Then  the set of distributions of $x^{k }_{\cdot}$ on 
$C([0,\infty),\bR^{d})$ is tight and any weakly converging
subsequence of distributions converges weakly to 
the distribution of one of solutions of
\eqref{11.29.2} as described in (i).
\end{theorem}

This theorem is proved in Section \ref{section 10.26.1}
by using Skorokhod's embedding method.  

Once the solvability of \eqref{11.29.2} is established,
the question of its weak uniqueness arises. 
A standard way (but there are other ways as well) of treating
it is the following when we argue formally without
caring about rigorousness at the moment.

Introduce
$$
L=(1/2)a^{ij}D_{ij}+b^{i}D_{i}.
$$

Take smooth $f$, $\lambda>0$, $R>|x |$, and find
a bounded sufficiently regular   solution of
\begin{equation}
                                     \label{10.18.1}
\lambda u-Lu=f
\end{equation}
in $B_{R}$ ($B_{\infty}:=\bR^{d}$) with zero boundary data (no boundary data 
if $R=\infty$). For $n=0,1,...$, let   
$g_{n}(y_{0},...,y_{n})$, $y_{k}\in\bR^{d}$, $k=0,1,...,n$,
be smooth bounded functions. Let
  $0=t_{0}\leq t_{1}<t_{2}...<\infty$.
 Use It\^o's formula to get that on the set $\{t_{n}\leq\tau\}$,
where $\tau=\tau_{R}$ is the first exit time of $x_{t}$ from $B_{R}$,
$$
u(x_{t\wedge\tau })e^{-t\wedge\tau }=u(x_{t_{n}})e^{-\lambda t_{n}}-
\int_{t_{n}}^{t\wedge\tau }e^{-\lambda s}f(x_{s})\,ds
$$
$$
+\int_{t_{n}}^{t\wedge\tau}e^{-s}Du(x_{s})\sqrt{ a(x_{s})}\,dw_{s}.
$$
Take expectations of both sides multiplied by the indicator
of $\{t_{n}<\tau \}$ times
$g_{n}(x_{t_{0}},...,x_{t_{n}})$. The expectation containing the
stochastic integral, naturally, disappears and after letting $t
\to\infty$ we obtain
\begin{equation}
                                                    \label{10.20.1}
Eg_{n}(x_{t_{0}},...,x_{t_{n}}) u(x_{t_{n}})I_{t_{n}<\tau }
e^{-\lambda t_{n}}
=\int_{t_{n}}^{\infty}e^{-\lambda s}Eg_{n}(x_{t_{0}},...,x_{t_{n}}) 
 f(x_{s})I_{s<\tau }\,ds.
\end{equation}
If the left-hand side is uniquely defined
(that is, independent of which solutions $x_{t}$ we take)
 for any $\lambda>0$ and smooth $f$, then, 
for $s\geq t_{n}$ and smooth $f$,
\begin{equation}
                                                    \label{10.20.2}
Eg_{n}(x_{t_{0}},...,x_{t_{n}}) 
 f(x_{s})I_{s<\tau_{R}}\quad\text{and hence}\quad
Eg_{n+1}(x_{t_{0}},...,x_{t_{n+1}}) I_{t_{n+1}<\tau_{R}}
\end{equation}
are uniquely defined. For $n=0$ the left-hand side of
\eqref{10.20.1} is $Eg_{0}u(x _{0})=g_{0}u(x)$ that is independent
of which solution
we take. By induction this allows us to conclude that
all quantities in \eqref{10.20.2} are uniquely defined
for all $n$. Letting $R\to\infty$
(if the above $R<\infty$) yields weak uniqueness.

There are the following obstacles to implementation of this scheme
if $b\in L_{d}$.
  To apply  It\^o's formula we,
generally, need $u\in W^{2}_{d}$
and $Lu \in L_{d}$. Then $b^{i}D_{i}u$ should be
in $L_{d}$, but $b$ is only in $L_{d}$ and then, apparently,
$Du$ needs to be bounded. However, by
embedding theorems the fact that $u\in W^{2}_{d}$
does not imply that $Du$  is bounded.
On the other hand, if  $u\in W^{2}_{p}$ with $p<d$
and $b\in L_{d}$ then $b^{i}D_{i}u\in L_{p}$
and the above mismatch does not occur, but
we need to know that It\^o's formula
is applicable to $u\in W^{2}_{p}$ for some $p<d$.

However, if $b$ is bounded, It\^o's formula is applicable
to $u\in W^{2}_{d-\varepsilon}$ for some $\varepsilon>0$
 (a consequence of a Fabes-Stroock result from
\cite{FS_84} if $b\equiv0$, which was carried over to bounded $b$
 in \cite{Ca_95} and to $b\in L_{d+\varepsilon}$
in \cite{Fo_98}). We would be in business
if we knew that this also holds if $b\in L_{d}$.
Then   weak uniqueness would follow from the solvability
of \eqref{10.18.1} in $W^{2}_{d-\varepsilon}$   for small
$\varepsilon>0$.

For $R\in(0,\infty)$ introduce
$\WO^{2}_{p}(B_{R})=W^{2}_{p}(B_{R})\cap
\{u:u|_{\partial B_{R}}=0\}$, $\WO^{2}_{p}(B_{\infty})
=W^{2}_{p}(\bR^{d})$,
$$
L_{0}=(1/2)a^{ij}D_{ij}.
$$
In a subsequent article the author intends to prove the following.
\begin{theorem}
                                     \label{theorem 10.18.2}
There is a $d_{0}=d_{0}(d,\delta,\|b\|) \in(d/2,d)$ such that, if $p\in[d_{0},d)$
and  for some
$R, \lambda>0$, and 
any $t\in[0,1]$ and $u\in \WO^{2}_{p}(B_{R})$ we have
\begin{equation}
                                     \label{10.18.3}
\|u\|_{W^{2}_{p}(B_{R})}\leq N\|\lambda u
-(t\Delta+(1-t)L_{0})u\|_{L_{p}(B_{R})}
\end{equation}
with $N$ independent of $u$ and $t$, then for any smooth $f$
equation \eqref{10.18.1}  has a unique solution
in $\WO^{2}_{p}(B_{R})$.
\end{theorem}

To the best of the author's knowledge the most general
conditions on the coefficients $a^{ij}$ when \eqref{10.18.3} holds
with $R=\infty$
(for  $p>2$) is given in \cite{Kr_09}, where
the solvability in $W^{2}_{p} $ spaces
is proved for second-order elliptic
equations with coefficients which are measurable
in one direction and VMO in the orthogonal directions
in each small ball with the direction depending on the ball.
Of course, we know from \cite{ADN_59} that 
\eqref{10.18.3} holds for all $p>1$ if $a^{ij}$
are continuous. In that case weak uniqueness with bounded $b$ 
was first proved by
 Stroock and Varadhan (see \cite{SV_79}).
If $a^{ij}\in VMO$, \eqref{10.18.3} for all $p>1$ and $R<\infty$
is proved in \cite{CFL_93}.

Most likely the conclusion of Theorem \ref{theorem 10.18.2}
is false if $p\geq d$   even if $a^{ij}=\delta^{ij}$. If $p<d$,
the proof of Theorem \ref{theorem 10.18.2} uses the fact that
$$
\int_{\bR^{d}} b^{p}|Du|^{p}\,dx\leq\Big
(\int_{\bR^{d}} b^{d}\,dx\Big)^{p/d}\Big(\int_{\bR^{d}} 
|Du|^{pd/(d-p)}\,dx\Big)^{(d-p)/d}
$$
$$
\leq N\int_{\bR^{d}} |D^{2}u|^{p}\,dx,
$$
where the last inequality follows from the embedding theorem
saying that $Du\in L_{q}$ if $D^{2}u\in L_{p}$ and
$$
2-\frac{d}{p}=1-\frac{d}{q}, 
\quad q=\frac{pd}{d-p}.
$$

It is worth saying that
Ladyzhenskaya-Ural'tseva in \cite{LU_73}
studied the case of $b\in L_{d+\varepsilon}$ with $p=2$
in \eqref{10.18.3} and continuous $a^{ij}$.
Actually in this situation
the assumption of Theorem \ref{theorem 10.18.2} is satisfied
 for any $p\in(1,\infty)$ 
 and hence the above described method of proving
  weak uniqueness works. But the case $b\in L_{d}$
is excluded.
 In the classical book
Gilbarg-Trudinger \cite{GT_01}
integrable drifts are not treated.

The last ingredient in the above scheme of how to prove
weak uniqueness on the basis of Theorem  \ref{theorem 10.18.2} is
 It\^o's formula. To state it we need some notation 
and   
{\em assumptions used throughout the paper\/}.

Let  $d_{1}$ be an   integer $>1$,
$(\Omega,\cF,P)$ be a complete probability space,
and let $(w_{t},\cF_{t})$ be a $d_{1}$-dimensional
Wiener process on this space with complete, relative to
$\cF,P$, $\sigma$-fields $\cF_{t}$. Let $\sigma_{t},t\geq0$,
be a progressively measurable process with values in the set 
of $d\times d_{1}$-matrices and let $b_{t},t\geq0$, be an $\bR^{d}$-valued
progressively  measurable process.
Assume that for any $T\in[0,\infty)$ and $\omega$
\begin{equation}
                                             \label{8.19.2}
\int_{0}^{T}\big(|\sigma_{t}|^{2}+|b_{t}|)\,dt<\infty.
\end{equation}

Under this condition the stochastic process  
\begin{equation}
                                             \label{8.19.1}
x_{t}=\int_{0}^{t}\sigma_{s}\,dw_{s}
+\int_{0}^{t}b_{s}\,ds
\end{equation}
is well defined. Fix a nonnegative Borel $\gb$ on $\bR^{d}$
and $\delta\in(0,1)$.

\begin{assumption}
                                   \label{assumption 11.15.1} 
We have $\|\gb\|:=\|\gb\|_{L_{d}(\bR^{d})}<\infty$ and
\begin{equation}
                                             \label{11.15.2}
|b_{t}|\leq \gb(x_{t}) ,\quad a _{t}\in \bS_{\delta}
\end{equation}
for  all $(\omega,t)$, where $a_{t}=\sigma_{t}
\sigma^{*}_{t}$.

\end{assumption}
 
Introduce
$$
L_{t}u(x_{t})=(1/2)a^{ij}_{t}D_{ij}(x_{t})+b^{i}_{t}D_{i}(x_{t}).  
$$

\begin{theorem}
                                         \label{theorem 9.6.2}
 
Under   Assumption  \ref{assumption 11.15.1}
there is a $d_{0}=d_{0}(d,\delta,\|\gb\|) \in(d/2,d)$ such that
 if    $p\geq d_{0}$ and
 $u\in W^{2}_{p}(\bR^{d})$,
then (a.s.) for all $t\geq0$
\begin{equation}
                                              \label{9.6.6}
u(x_{t})=u(0)+\int_{0}^{t}L_{s}u(x_{s})\,ds
+\int_{0}^{t}D_{i}u(x_{s})\sigma^{ik}_{s}\,dw^{k}_{s}
\end{equation}
and the last term is a square integrable martingale. 
\end{theorem}
This theorem is proved in Section \ref{section 10.26.1}.

The above results and the discussion after Theorem
\ref{theorem 9.6.4} immediately yield the following.

\begin{theorem}
                                         \label{theorem 10.27.1}
In the setting of Theorem \ref{theorem 9.6.4} suppose
that for a  $p\in[d_{0},d)$
and  any
$R, \lambda\in(0,\infty)$,   
 $t\in[0,1]$, $u\in \WO^{2}_{p}(B_{R})$ we have
$$
\|u\|_{W^{2}_{p}(B_{R})}\leq N\|\lambda u
-(t\Delta+(1-t)L_{0})u\|_{L_{p}(B_{R})}
$$
with $N$ independent of $u$ and $t$ (may depend on $\lambda,R$,...).
Then solutions of \eqref{11.29.2} are weakly unique.
\end{theorem}

In the heart of the above results lies the following
estimate.

\begin{theorem}
                                           \label{theorem 10.19.2}
Under  Assumption  \ref{assumption 11.15.1}
there is a $d_{0}  \in(d/2,d)$, depending
only on $d,\delta,\|\gb\|$, such that,
for any $\lambda>0$, $p\geq d_{0}$,
and nonnegative Borel $f(x)$ given on $\bR^{d}$  we have
\begin{equation}
                                          \label{9.5.6}
E\int_{0}^{\infty}e^{-\lambda t}\Psi_{\lambda}(x_{t})
f(x_{t})\,dt\leq N \lambda^{d/(2p)-1}\|f\|_{L_{p}(\bR^{d})}.
\end{equation}
where $N$ depends only on $d,\delta$, and $\|\gb\|$,
$\Psi_{\lambda}(x)=\exp( \sqrt{\lambda}\nu|x|)$, $\nu=\mu/4$,
and $\mu$ is taken from
 Theorem \ref{theorem 9.3.1}. 
\end{theorem}

This theorem is proved in Corollary \ref{corollary 9.5.1}.

The above results allows one to construct
Markov diffusion processes corresponding to $L$.
To show how to do this we need the following,
which would be a simple consequence of Theorem 4.5.1
of \cite{SV_79} were $b$   supposed to be bounded.

\begin{lemma}
                                            \label{lemma 12.6.1}
Let $a$ and $b$ be the same as before Theorem \ref{theorem 9.6.4}.
Suppose that we are given a continuous process
$x_{t}$, $t\geq0$, such that $x_{0}=0$, for any $T\in(0,\infty)$
$$
\int_{0}^{T}|b(x_{t})|\,dt<\infty
$$
(a.s.), and for any
twice continuously differentiable function $u(x)$
with compact support the process
\begin{equation}
                                               \label{11.6.2}
u(x_{t})-\int_{0}^{t}Lu(x_{s})\,ds
\end{equation}
is a local martingale with respect  to the filtration
of $\sigma$-fields $\{\cN_{t}=\sigma(x_{s};s\leq t), t\geq0\}$.
Then there exists a $d$-dimensional Wiener process
$(w_{t},\bar \cN_{t})$, $t\geq0$, where $\bar \cN_{t}$ is the completion
of $\cN_{t}$, such that \eqref{11.29.2} is satisfied with $x=0$.
\end{lemma}

Proof. First observe that by using cut-off functions
one easily shows that \eqref{11.6.2} is a  local martingale
for any twice continuously differentiable function $u$.
Then,
we claim that the following processes are local martingales
$$
X_{t}:=x_{t}-\int_{0}^{t}b(x_{s})\,ds,
$$
$$
B_{t}:=x_{t}x^{*}_{t}-\int_{0}^{t}\big(a(x_{s})+b(x_{s})x^{*}_{s}
+x_{s}b^{*}(x_{s})\big)\,ds,
$$
$$
A_{t}:=X_{t}X_{t}^{*}-\int_{0}^{t} a(x_{s})\,ds.
$$

Indeed, the first two processes are obtained from 
\eqref{11.6.2} for $u=x,xx^{*}$. Concerning the last one
introduce $\gamma_{R}$ as the minimum of $\tau_{R}$   
and 
$$
\inf\{t\geq0:\int_{0}^{t}|b(x_{s})|\,ds+|B_{t}|\geq R\}.
$$
Also let
$$
\Phi_{t}=\int_{0}^{t}b(x_{s})I_{s<\gamma_{R}}\,ds.
$$
Observe that $X_{t\wedge\gamma_{R}}$ and $\Phi_{t}$
are bounded and simple manipulations show that
$$
A_{t\wedge\gamma_{R}}=\int_{0}^{t}X_{s\wedge\gamma_{R}}\,d
\Phi^{*}_{s}-X_{t\wedge\gamma_{R}} 
\Phi^{*}_{t}+
\int_{0}^{t}\big(d
\Phi_{s}\big)X^{*}_{s\wedge\gamma_{R}}- 
\Phi_{t}X^{*}_{t\wedge\gamma_{R}}+B_{t\wedge\gamma_{R}},
$$
which by the Lemma from Appendix 2 of \cite{Kr_77}
shows that $A_{t\wedge\gamma_{R}}$ is a martingale.

By the above claim the quadratic variation
process of the local martingale $X_{t}$
is
$$
\int_{0}^{t}a(x_{s})\,ds.
$$
After that our assertion   
follows directly from
  Theorem III.10.8 of \cite{Kr_95}. The lemma is proved. \qed  

\begin{remark}
We used a result  from \cite{Kr_95},
where the initial definition (see there Definition II.8.2)
 of a martingale is different from commonly used
 and, owing to Doob's optional stopping theorem,
seemingly admits wider class of processes than the martingales
in the usual sense. However, just considering
stopping times taking only two values, one easily sees
that, actually, martingales from \cite{Kr_95} are martingales
in the classical sense.
 \end{remark}

\begin{theorem}
                                                 \label{theorem 12.6.1}
Let $a$ and $b$ be as in Lemma \ref{lemma 12.6.1}.
Then there exists a
continuous strong Markov process $X=(x_{t},\infty,\cM_{t},P_{x})$
(the terminology taken from \cite{Dy_63})
in $\bR^{d}$ such that 
  for any $x\in\bR^{d}$ and $t\geq0$
\begin{equation}
                                               \label{11.4.3}
E_{x} \int_{0}^{t}|b(x_{s})|\,ds<\infty 
\end{equation}
and for any
twice continuously differentiable function $u$
with compact support the process
\eqref{11.6.2} is a local martingale
relative to $P_{x}$.
Furthermore, $(x_{t},\infty,\cM_{t+},P_{x})$
is a Markov process.
\end{theorem}

Proof. The proof of this theorem follows
the proof of Theorem 3 of \cite{Kr_73}
and we only point out the necessary changes
related to the fact that, unlike \cite{Kr_73}
where $b$ is Borel bounded, our $b\in L_{d}$.

As in \cite{Kr_73} we set $\Omega=C([0,\infty),\bR^{d})$
and for $\omega=\omega_{\cdot}\in \Omega$ define
$x_{t}(\omega)=\omega_{t}$. Also set
$\cM_{t}=\cN_{t}=\sigma(x_{s};s\leq t)$ and by $\Pi_{x}$
denote the set of probability measures $P$ on
$(\Omega,\cN_{\infty})$ such that $P(x_{0}=x)=1$
and the process \eqref{11.6.2} is a local martingale
for any
twice continuously differentiable function $u$.
According to Theorem \ref{theorem 9.6.4}
and Theorem \ref{theorem 10.19.2}, assuring that
\eqref{11.4.3} holds for solutions of
\eqref{11.29.2}, $\Pi_{x}\ne\emptyset$.

Owing to Lemma \ref{lemma 12.6.1}, Corollary 1.2
of \cite{Kr_19} and
Theorem \ref{theorem 10.19.2} are applicable, that is,
for $P\in\Pi_{x}$ and
any $n\geq0$ 
\begin{equation}
                                             \label{10.20.10} 
 E\max_{r\in[s,t]}|x_{r}-x_{s}|^{2n}\leq N (t-s)^{n},
\end{equation}
where $N=N(n,d,\delta,\|b\|)$ and \eqref{9.5.6}
holds with $N=N(d,\delta,\| b\|)$. In particular,  
the assertions of
Lemmas 5 and 6 of \cite{Kr_73} are valid.
After that the proof goes the same way as in \cite{Kr_73}
once more using Theorem \ref{theorem 9.6.4},
this time its second statement, while proving that
$\{\Pi_{x}\}$ is a $B$-system in the terminology of \cite{Kr_73}.
The theorem is proved. \qed

\begin{remark}
                                                \label{remark 12.6.3}
In a subsequent article we will show that 
the process $(x_{t},\infty$, $\bar\cM_{t+},P_{x})$ is  strong Markov
with strong Feller semigroup (see \cite{Kr_20}).
\end{remark}

Theorem \ref{theorem 12.6.1} provides existence of 
a Markov diffusion process corresponding to the operator $L$.
One knows that, generally, this process is not unique
in any sense. In this connection we   present some results 
such as   Corollary \ref{corollary 10.11.1}
which are the main tools in proving the Harnack inequality
and H\"older continuity property of harmonic functions
for the corresponding diffusion processes with drift in $L_{d}$.

We also deal with some issues from the theory
of partial differential equations. For instance,
Corollary \ref{corollary 10.7.1}, in particular,
provides the maximum principle for elliptic equations with
measurable $a^{ij}$ and drift in $L_{d}$ for solutions
in $W^{2}_{p}$ with $p<d$ (in case $p=d$ this is a classical
Aleksandrov's result and, if $b$ is bounded, $p<d$ is allowed
according to the results in \cite{FS_84}, \cite{Fo_98}).
 Theorem \ref{theorem 10.14.2}
is indispensable
  in the Sobolev space 
theory of fully nonlinear elliptic
 equations while studying the possibility to pass  to the limit
in such equations.

The rest of the article is organized as follows.
In   Section \ref{section 10.19.3}    we prove 
Theorem \ref{theorem 10.19.2}. This allows us to prove Theorems
\ref{theorem 9.6.4} and \ref{theorem 9.6.2}
in Section \ref{section 10.26.1}. 
Section \ref{section 10.15.1} is devoted
to studying fine properties of our processes such as
estimating the time spent in sets of small measure,
the probability to reach such sets, 
Fanghua Lin estimates playing a major role in the Sobolev space theory
of {\em fully nonlinear elliptic equations\/},
boundary behavior of solutions of the corresponding
elliptic equations with first order coefficients in 
$L_{d}$, and the probability to pass through narrow tubes,
which in the first draft of the paper was one of cornerstones
of everything else. We also prove the doubling property
of the corresponding  Green's measures 
and the fact that their densities are in the class $A_{\infty}$
of Muckenhoupt.

We finish the introduction with some notation and the stipulation
about constants. In the proofs of various results  we use
the symbol $N$  to denote finite 
nonnegative constants
which may change from one occurrence to another and
we do not always specify on which data these  constants
depend. In these cases the reader should remember
that, if in the statement of a result there are constants
called $N$ which are claimed to depend only on certain
parameters, then in the proof of the result
the constants $N$ also depend only on the same
parameters unless specifically stated otherwise.
Of course, if we write 
$
N=N(...),
$
 this means that $N$ depends only
on what is inside the parentheses. 
Another point is that when we say that certain constants
depend only on such and such parameters we mean, in particular,
  that the dependence is such that these constants stay bounded
as the parameters vary in compact subsets of their ranges.

Introduce $|\Gamma|$ as the volume of $\Gamma\subset\bR^{d}$,
$$
\dashint_{B}f(x)\,dx=\frac{1}{|B|}\int_{B}f(x)\,dx,\quad 
a_{\pm}=a^{\pm}=(1/2)(|a|\pm a).
$$
Use  the notation $u^{(\varepsilon)}=u*\zeta_{\varepsilon}$,
where
$\zeta_{\varepsilon}(x)=\varepsilon^{-d}\zeta(x/\varepsilon)$,
$\varepsilon>0$,
and $\zeta$ is a nonnegative $C^{\infty}$-function with
support in $B_{1}$ whose integral is equal to one.
If $B$ is a ball and $\eta$ is a positive number, by $\eta B$
we denote a concentric ball whose radius is $\eta$
times that of $B$. 

If it is not explicitly stated otherwise, by $x_{t}$
we always mean the process defined by \eqref{8.19.1}
and let
$\tau_{R}(x)$ be the first exit time of $x+x_{t}$ from  
$B_{R}$ (equal to infinity if $x+x_{t}$ never exits from $B_{ R}$).
 Also let $\tau_{R}=\tau_{R}(0)$.
 
\mysection{Green's functions}

                                             \label{section 10.19.3}
 
We suppose throughout the article
 that Assumption \ref{assumption 11.15.1} 
is satisfied.
Recall that Theorem 2.17 of \cite{Kr_19} implies that
if $p\geq d$, 
 then there exists   constants $N$ and $\mu>0$,
depending only on $d,p,\delta$, and $\|\gb\|$, such that 
for any  
$\lambda>0$ and Borel nonnegative $f$ given on $\bR^{d}$
we have
\begin{equation}
                                                   \label{7.29.2}
 E\int_{0}^{\infty}e^{-\lambda t}
f(x_{t}) \,dt   \le N\lambda^{d/(2p)-1}
\|\Psi_{\lambda}^{-1} f\| _{L_{p}(\bR^{d})},
\end{equation}
where $\Psi_{\lambda}(x)=\exp( \sqrt{\lambda}\mu|x|)$.

Here is a straightforward consequence
of \eqref{7.29.2}.

\begin{theorem}
                                       \label{theorem 9.3.1}
Let   $p\geq d$.
 Then there exists   constants $N$ and $\mu>0$,
depending only on $d,\delta,p$, and $\|\gb\|$, and for any
$\lambda>0$
  there exists a nonnegative Borel function $G_{\lambda}(x)$
(Green's function of $x_{\cdot}$)
on $\bR^{d}$
 such that  for any Borel nonnegative $f$ given on $\bR^{d}$
we have
$$
 E\int_{0}^{\infty}e^{-\lambda t}
f(x_{t}) \,dt =\int_{\bR^{d}}f(x)G_{\lambda}(x)\,dx,
$$
\begin{equation}
                                                   \label{9.3.5}
 \|\Psi_{\lambda} G_{\lambda}\|_{L_{p/(p-1)}(\bR^{d})} 
\le N\lambda^{d/(2p)-1}.
\end{equation}
 
\end{theorem}

The highest power of summability of $G_{\lambda}$ guaranteed by this theorem
occurs when $p=d$ and this is $d/(d-1)$. It turns out that,
actually,
$G_{\lambda}$ is summable to a higher power. The proof of this
is based on Gehring's lemma from \cite{Gh_73},
Aleksandrov's estimates, and the following,
which is a particular case of  Lemma 2.13 in \cite{Kr_19}

\begin{lemma}
                                             \label{lemma 9.2.1}
There is a constant $N=N(d,\delta,\|\gb\|)$ such that
for any $R\in(0,\infty)$
\begin{equation}
                                          \label{9.2.3}
NE\int_{0}^{\tau_{R}\wedge R^{2}}e^{-t}\,dt\geq R^{2}\wedge 1.
\end{equation}
 
\end{lemma}

\begin{theorem}
                                           \label{theorem 9.3.2}
There are constants $\varepsilon\in(0,1)$ and $N $, depending only
on $d,\delta$, and $\|\gb\|$, such that for any   ball $B$
of radius $R\leq 1/2$ and $p\geq d_{0}: =d-\varepsilon$, we have
\begin{equation}
                                          \label{10.14.1}
\| G_{1}\|_{L_{p/(p-1)}(B)}\leq N R^{-d /p }
\| G_{1}\|_{L_{1}(2B )} ,
\end{equation}
which is equivalently rewritten as
$$
\Big(\dashint_{B}G^{p/(p-1)}_{1}\,dx\Big)^{(p-1)/p}
\leq N\dashint_{2B }G_{1}\,dx.
$$
 
\end{theorem}

Proof. We basically follow the arguments in \cite{FS_84}.
Take $R\in(0,1/2]$, a closed ball $B$ of radius $R$
and let $B'$ be the concentric open ball of radius $2R$.
Define 
recursively 
$$
\gamma^{1}=\inf\{t\geq0 :x_{t}\in B\},\quad
\tau^{1}=\inf\{t\geq\gamma^{1} :x_{t}\not\in B'\},
$$
$$
\gamma^{n+1}=\inf\{t\geq\tau^{n} :x_{t}\in B\},\quad
\tau^{n+1}=\inf\{t\geq\gamma^{n+1} :x_{t}\not\in B'\}.
$$
Then for any nonnegative Borel $f$ vanishing outside $B$
with $\|f\|_{L_{d}(B)}=1$
we have
$$ 
\int_{B}f   G_{1}(x)\,dx
=E\int_{0}^{\infty}e^{-t}f (x_{t})\,dt
$$
$$
=\sum_{n=1}^{\infty}Ee^{-\gamma^{n}}
E\Big(\int_{\gamma^{n}}^{\tau^{n}}e^{-(t-\gamma^{n})}f(x_{t})\,dt
\mid \cF_{\gamma^{n}}\Big).
$$
Next we use the conditional  version of the Aleksandrov
estimate to see that the conditional expectation
above is less than $NR\|f\|_{L_{d}(B')}=NR$. After that
we use the conditional  version of \eqref{9.2.3}
to get that
$$
R^{2}\leq N E\Big(\int_{\gamma^{n}}^{\tau^{n}}e^{-(t-\gamma^{n})} \,dt
\mid \cF_{\gamma^{n}}\Big).
$$
Then we obtain
$$
\int_{B}f   G_{1}(x)\,dx
\leq N  R ^{-1}
\sum_{n=1}^{\infty}E e^{-\gamma^{n}}
E\Big(\int_{\gamma^{n}}^{\tau^{n}}e^{-(t-\gamma^{n})} \,dt
\mid \cF_{\gamma^{n}}\Big)
$$
$$
=N  R ^{-1}
\sum_{n=1}^{\infty}E  
 \int_{\gamma^{n}}^{\tau^{n}}e^{- t } \,dt
$$
$$
\leq N  R ^{-1}
E \int_{0}^{\infty}e^{-t}I_{B'}(x_{t}) \,dt
=N  R ^{-1}\int_{B'}G_{1}(x)\,dx.
$$

The arbitrariness of $f$ 
implies that
$$
\Big(\dashint_{B}  G_{1} ^{d/(d-1)}(x)\,dx
\Big)^{(d-1)/d}\leq N \dashint_{B'}  G_{1}(x)\,dx.
$$

Now the assertion of the theorem for $p=d_{0}$ follows directly from
the corrected version of the famous Gehring's lemma
proved as Proposition 5.1  in \cite{GM_79}. 
For larger $p$ it suffices to use H\"older's inequality.
The theorem is proved.  \qed

\begin{remark}
                                           \label{remark 10.14.1}
Take $\nu\in\bR$. Then the assertion of Theorem
\ref{theorem 9.3.2} remains true if we replace
\eqref{10.14.1} with
\begin{equation}
                                          \label{9.3.05}
\|\Psi G_{1}\|_{L_{p/(p-1)}(B)}\leq N R^{-d /p }
\|\Psi G_{1}\|_{L_{1}(2B )} ,
\end{equation}
where $\Psi(x)=\exp(\nu|x|)$, and allow $N$ to depend also  on $\nu$.
This follows from the fact that the supremum of $\Psi$
over $B$ is less than a constant independent of $R\leq1/2$ times
the infimum of $\Psi$ over $2B $.
\end{remark}

Here is a substantial improvement of Theorem
\ref{theorem 9.3.1}. Below and many times
in the future  we use
self-similarity transformations like $x_{t}\to
cx_{t/c^{2}}$, where $c>0$ is a constant. This 
transformation changes $\sigma_{t}$ and $b_{t}$
in a well known way, which will bring about a new
function $\gb$ (see \eqref{11.15.2}). A remarkable
fact is that this new $\gb$ has {\em the same\/}
$L_{d}(\bR^{d})$-norm as the original one.

\begin{corollary}
                                           \label{corollary 9.5.1}
For any $\lambda>0$ and $p\geq d_{0}$ we have
\begin{equation}
                                                   \label{9.5.2}
 \|\Psi _{\lambda}G_{\lambda}\|_{L_{ p/(p-1) }(\bR^{d})}
\le N\lambda^{d/(2p)-1},
\end{equation}
where $N$ depends only on $p,d,\delta$, and $\|\gb\|$,
$\Psi_{\lambda}(x)=\exp( \sqrt{\lambda}\nu|x|)$, $\nu=\mu/4$,
and $\mu$ is taken from
 Theorem \ref{theorem 9.3.1}. In other words,
for any nonnegative Borel $f(x)$ given on $\bR^{d}$
estimate \eqref{9.5.6} holds:
$$
E\int_{0}^{\infty}e^{-\lambda t}\Psi_{\lambda}(x_{t})
f(x_{t})\,dt\leq N \lambda^{d/(2p)-1}\|f\|_{L_{p}(\bR^{d})}.
$$
\end{corollary}

Indeed, the case of arbitrary $\lambda>0$ reduces
to the one where $\lambda=1$ by using 
self-similarity
 and for $\lambda=1$
it suffices to note that, for  $q= p/(p-1)$,
$$
\int_{\bR^{d}}e^{q\nu |x|}G_{1}^{q}(x)\,dx=N
\int_{\bR^{d}}\Big[\int_{|y-x|\leq1}e^{q\nu | y|}G_{1}^{q}( y)
\,dy\Big]\,dx
$$
$$
\leq N\int_{\bR^{d}}e^{-q\nu|x|}
\Big[\int_{|y-x|\leq1}e^{q2\nu | y|}G_{1}^{q}( y)
\,dy\Big]\,dx
$$
and then, to estimate the interior
integral, use \eqref{9.3.05} and the fact that
owing to \eqref{9.3.5}
and H\"older's inequality
\begin{equation}
                                                   \label{9.5.1}
 \int_{\bR^{d}}e^{ 2\nu | x|}  G_{1}(x)\,dx 
\le N.
\end{equation}

\begin{corollary} 
                                            \label{corollary 9.7.2}
For any $x\in\bR^{d}$, $p\geq d_{0}$, 
and Borel nonnegative $f$ vanishing outside
  $B_{1}(x)$   we have
$$
E\int_{0}^{\infty}e^{-t}f( x_{t})\,dt 
  \le N e^{-\nu|x|}
\|f\| _{L_{p}(B_{1}(x))},
$$
where $N$   depends  only on $p,d,\delta$, and $\|\gb\|$ 
and $\nu$ is the same as in Corollary \ref{corollary 9.5.1}.

\end{corollary}

\begin{theorem}
                                          \label{theorem 9.7.1}
There is a constant $N=N(p,d,\delta,\|\gb\|)$ such that for any
$n=1,2,...$, nonnegative Borel $f$ on $\bR^{d}$,
 $T\in(0,\infty)$,  and $p\geq  d_{0}$  we have
\begin{equation}
                                          \label{9.7.1}
E\Big[\int_{0}^{T}  
f(x_{t})\,dt\Big]^{n}\leq n!N^{n} 
T^{n[1-d/(2p)]}\|\Psi^{- 1/n }_{1/T}f\|^{n}_{L_{p}(\bR^{d})},
\end{equation}
where $\Psi_{\lambda}$ is taken from Corollary
\ref{corollary 9.5.1}.
\end{theorem}

Proof. We are going to use the induction on $n$. 
Our induction hypothesis is that for an $n=1,2,...$, any
nonnegative Borel $f$,   $x\in\bR^{d}$,
and $\kappa\in[0,1/n]$
  \begin{equation}
                                        \label{8.30.1}
E\Big(\int_{0}^{ T}f(x+x_{t})\,dt\Big)^{n}\leq n!
N^{n}T^{n(1-d/(2p))}\Psi_{1/T}^{\kappa n}(x)\|\Psi^{-\kappa}_{1/T}f\|^{n}
_{L_{p}(\bR^{d})}
\end{equation}

If the hypothesis  holds true for some $n\geq1$,
then by using its conditional  version
and the fact that
$$
I:=E\Big(\int_{0}^{T} 
f(x+x_{t})\,dt\Big)^{n+1} 
$$
$$
=(n+1)E\int_{0}^{T}f(x+x_{t})E\Big\{\Big[\int_{0}^{T-t} 
f(x+x_{t}+(x_{t+r}-x_{t}))\,dr\Big]^{n}\mid\cF_{t}\Big\}\,dt,
$$
we see that, for any $\kappa\in(0,1/n)$,
\begin{equation}
                                        \label{8.31.1}
I\leq (n+1)!N^{n}T^{n[1-d/(2p)]}\|\Psi^{-\kappa}_{1/T}
f\|^{n}_{L_{p}(\bR^{d})}
E\int_{0}^{T}\Psi^{\kappa n}_{1/T}f(x+x_{t})\,dt.
\end{equation}

Next, introduce
$$
F(T)=E\int_{0}^{T}\Psi^{\kappa n}_{1/T}f(x+x_{t})\,dt
$$
 and observe that for any  $\lambda>0$
owing to Corollary \ref{corollary 9.5.1} we have
$$
  e^{-\lambda T}F(T)
\leq \lambda\int_{0}^{\infty}F(t)e^{-\lambda t}\,dt
$$
$$
=E\int_{0}^{\infty}e^{-\lambda t} \Psi^{\kappa n}_{1/T}
f(x+x_{t})\,dt\leq N \lambda^{d/(2p)-1}\| 
\Psi^{-\mu}_{\lambda}(\cdot-x)\Psi^{\kappa n}_{1/T}f
\|_{L_{p}(\bR^{d})}
$$
$$
\leq N \lambda^{d/(2p)-1}\Psi^{\mu} _{\lambda}( x)\| 
 \Psi^{-\mu} _{\lambda}\Psi^{\kappa n }_{1/T}f
\|_{L_{p}(\bR^{d})},
$$
where $\mu\in[0,1]$ and the last inequality is due to the fact that
$\Psi^{-1}_{\lambda}(y-x)\leq \Psi^{-1} _{\lambda}(y)
\Psi _{\lambda}( x)$.
For $\lambda=1/T$ we get \eqref{8.30.1} with $n=1$,
which justifies the start of the induction.
For $\mu=\kappa(n+1)$, $\kappa\in[0,1/(n+1)]$, we have
$\Psi^{-\mu} _{1/T}\Psi^{\kappa n }_{1/T}=
\Psi^{-\kappa   }_{1/T}$ and this along with 
\eqref{8.31.1} show  that our hypothesis holds true
also for $n+1$.

Now it only remains to observe that for $x=0$
and $\kappa=1/n$ estimate \eqref{8.30.1} coincides
with \eqref{9.7.1}. 
  The theorem is 
proved. \qed

Next theorem improves Theorem 1.1 of \cite{Kr_19} 
in what concerns the range of $p$
for uniformly nondegenerate processes.

\begin{theorem}
                                       \label{theorem 9.5.1}

There exists a constant $N=N(p,d,\delta,\|\gb\|)$ such that,
for any $R>0$,   $p\geq d_{0}$, and Borel nonnegative $f$ 
given on $B_{R}$,
we have
\begin{equation}
                                                     \label{9.5.4}
E\int_{0}^{\tau_{R}}f(x_{t})\,dt\leq
NR^{2-d/p}\|f\|_{L_{p}(B_{R})}.
\end{equation}

\end{theorem}

Proof. H\"older's inequality allows us to only concentrate
on $p=d_{0}$.
Scalings show that we may assume that $R=1$.
Also we may assume that $f$ is bounded and is zero outside $B_{1}$.
In that case denote by $\frM$ the set of stopping times
$\gamma\leq\tau:=\tau_{1}$, and set
$$
u_{\gamma}=E\Big[\int_{\gamma}^{\tau}f(x_{t})\,dt
\mid \cF_{\gamma}\Big],\quad  \bar u=\esssup_{\gamma\in\frM}u_{\gamma}.
$$
Observe that for any $\omega$  and $\lambda>0$
it holds that
$$
 \int_{\gamma}^{\tau}f(x_{t})\,dt=
\int_{\gamma}^{\tau}e^{-\lambda (t-\gamma)}f(x_{t})\,dt
+\lambda \int_{0}^{\infty}e^{-\lambda (t-\gamma)}
I_{\gamma\leq t<\tau}\Big[\int_{t}^{\tau}f(x_{s})\,ds\Big]\,dt.
$$
By the conditional  version of \eqref{9.5.6} 
(recall that $p=d_{0}$) (a.s.)
$$
E\Big[\int_{\gamma}^{\tau}e^{-\lambda (t-\gamma)}f(x_{t})\,dt
\mid \cF_{\gamma}\Big]
$$
$$
\leq
E\Big[\int_{\gamma}^{\infty}e^{-\lambda (t-\gamma)}f(x_{t})\,dt
\mid \cF_{\gamma}\Big]\leq 
N \lambda^{d/(2p)-1}\|f\|_{L_{p}(B_{1})}.
$$
Hence,
$$
u_{\gamma}\leq N \lambda^{d/(2p)-1}\|f\|_{L_{p}(B_{1})}
$$
$$
+\lambda E\Big[\int_{\gamma}^{\infty}e^{-\lambda (t-\gamma)}
I_{\gamma\leq t<\tau}E\Big[\int_{t}^{\tau}f(x_{s})\,ds
\mid \cF_{t}\Big]\,dt
\mid \cF_{\gamma}\Big],
$$
where the last term is dominated by
$$
\lambda \bar u E\Big[\int_{\gamma}^{\infty}e^{-\lambda (t-\gamma)}
I_{\gamma\leq t<\tau} \,dt
\mid \cF_{\gamma}\Big]\leq
\lambda \bar u E\Big[\int_{\gamma}^{\tau} \,dt
\mid \cF_{\gamma}\Big]\leq N_{1}\lambda \bar u
$$
(a.s.), where the last inequality follows from 
the  conditional  version of Corollary 2.1 of \cite{Kr_19}.
Thus, (a.s.)
$$
u_{\gamma}\leq N \lambda^{d/(2p)-1}\|f\|_{L_{p}(B_{1})}
+N_{1}\lambda \bar u.
$$
Since $\gamma$ is arbitrary within $\frM$,
\begin{equation}
                                                       \label{9.5.7}
\bar u \leq N \lambda^{d/(2p)-1}\|f\|_{L_{p}(B_{1})}
+N_{1}\lambda \bar u
\end{equation}
(a.s.), and since $\bar u<\infty$ ($f$ is bounded), by taking $\lambda
=1/(2N_{1})$, we arrive at 
$$
\bar u\leq 
N  \|f\|_{L_{p}(B_{1})}.
$$
The theorem is proved.\qed

\begin{remark}
                                            \label{remark 9.5.1}
Equation \eqref{9.5.7} implies
that  
$(p=)d_{0}\geq d/2$. Of course, the example of the Wiener
process with no drift shows more than that, namely, $d_{0}> d/2$.
\end{remark}

We finish the section with a result which will be used
in a subsequent article (see \cite{Kr_20}).

\begin{theorem}
                                        \label{theorem 11.23.1}
Let   $p\geq d_{0}$.
 Then there exists   constants $N$ and $\mu>0$,
depending only on $d,p$, and $\|\gb\|$,
and there exists $R_{0}=R_{0}(d,\|\gb\|)\geq2$, such that for any
$\lambda>0$, $R\in[0,\infty)$,
 and Borel nonnegative $f$ given on $\bR^{d}$
we have
\begin{equation}
                                                   \label{11.23.1}
E\int_{0}^{\infty}e^{-\lambda \phi_{t}(B_{R}^{c})}
f(x_{t}) \,dt   \le N (R\sqrt{\lambda}+R_{0}) ^{2-d/p} 
\lambda^{d/(2p)-1} 
\|\Psi_{R,\lambda}^{-1} f\| _{L_{p}(\bR^{d})},
\end{equation}
where $\Psi _{R,\lambda}(x)=\exp\big(\sqrt{\lambda}\,\mu\,
\dist(x,B_{R+R_{0}/\sqrt{\lambda}})\big)$ and
$$
\phi _{t}(B_{R}^{c})=\int_{0}^{t} I_{|x_{s}|\geq R}\,ds.
$$
\end{theorem}

This theorem looks very much like 
Theorem 2.18
of \cite{Kr_19} proved for possibly degenerate processes
for $p\geq d$ rather than $p\geq d_{0}$.
Theorem \ref{theorem 11.23.1}
is proved in the same way as Theorem 2.18
of \cite{Kr_19} on the basis of sharper estimates
of Green's functions in the special case of
uniformly nondegenerate processes. We only need to use again
Lemma 2.8 of \cite{Kr_19} and use our Theorem  \ref{theorem 9.5.1}
and Corollary \ref{corollary 9.5.1} instead of Theorems 1.1 and 2.17
of \cite{Kr_19}, respectively.

\mysection{It\^o's formula and solvability
of stochastic equations with drift in $L_{d}$}
                                         \label{section 10.26.1}

Recall that 
Assumption  \ref{assumption 11.15.1}
is supposed to be satisfied throughout the article.
First we deal with It\^o's formula.

{\bf 
Proof of Theorem \ref{theorem 9.6.2}}. Since $p>d/2$ (see Remark \ref{remark 9.5.1}),
by embedding theorems, $u$ is bounded and continuous.
 Furthermore, since $p\geq d_{0}$,   by embedding theorems,
the $L_{q}$-norms of $|Du|$ over any ball of radius
one are bounded by the same constant, where $q=d_{0}d/(d-d_{0})$.
Next since $2d_{0}\leq q$ ($d_{0}>d/2$),  
 for any $\lambda>0$ and $\Psi_{\lambda}$
from Corollary \ref{corollary 9.5.1}, it holds that
 $\Psi_{\lambda}^{-1}|Du|^{2}\in L_{d_{0}}(\bR^{d})$.
Therefore, by Theorem~\ref{theorem 9.7.1}
$$
E\int_{0}^{T}|D u(x_{t})|^{2}\,dt<\infty
$$
for any $T\in(0,\infty)$, which proves that the stochastic integral
in \eqref{9.6.6} is indeed a square integrable martingale.

We prove \eqref{9.6.6} by passing to the limit from smooth
functions $u_{n}$ which converge to $u$ in $W^{2}_{p}(\bR^{d})$.
In light of what is said in the previous paragraph,
$u_{n}\to u$ uniformly in $\bR^{d}$ and also there is no
difficulty to pass to the limit in the stochastic term.
In the deterministic term there could be only one  expression
of concern  
$$
E\int_{0}^{T}|b_{t}||D(u_{n}-u)|(x_{t})\,dt
$$
which owing to the condition $|b_{t}|\leq  \gb(x_{t})$
and Theorem \ref{theorem 9.7.1} is less than a constant
independent of $n$ times
$$
\|\Psi^{-1}_{1/T}\gb D(u_{n}-u) \,\|_{L_{d_{0}}(\bR^{d})}.
$$
The latter by H\"older's inequality is estimated by
$$
\|\gb\|\,\|\Psi^{-1}_{1/T}D(u_{n}-u) \|_{L_{q}(\bR^{d})},
$$
where the last term, by embedding theorem, is less than a constant
independent of $n$ times the $W^{2}_{p}(\bR^{d})$-norm of
$u_{n}-u$ that tends to zero as $n\to \infty$. This proves 
the theorem. \qed

Next we   deal with stochastic equations
with drift in $L_{d}$.

{\bf Proof of Theorem \ref{theorem 9.6.4}}. 
Having in mind mollifiers we see that assertion (ii)
implies (i).
The proof of (ii) is achieved by 
repeating the proof
of Theorem 2.6.1 of \cite{Kr_77} with only a few changes
which we point out below.  
By Corollary 1.2 of \cite{Kr_19}  
for any $m=1,2,...$, $0\leq s\leq t$
$$
E\max_{r\in[s,t]}|x^{(k)}_{r}-x^{(k)}_{s}|^{2m}\leq N (t-s)^{m},
$$
where $N$ is independent of $k$.  This yields
the tightness of distributions.

Then in exactly the same way as in the proof 
of Theorem 2.6.1 of \cite{Kr_77}, for any
weakly converging sequence $\{k'\}$ of distributions
of $x_{\cdot}^{(k')}$  by using Skorokhod's
embedding theorem, we find a probability space $(\Omega,\cF,P)$,
$d$-dimensional Wiener processes $(\tilde w^{(k')}_{t},\tilde
\cF^{(k')}_{t})$
defined on this space, and $\tilde\cF^{(k')}_{t}$-adapted
continuous processes $\tilde x^{(k')}_{t}$ such that,
for some $( \tilde  x _{t}, \tilde  w _{t})$ we have
  $( \tilde x^{(k')}_{t},\tilde w^{(k')}_{t})
\to ( \tilde  x _{t}, \tilde  w _{t})$ in probability
for any $t\geq0$ and for any $k'$ it holds that with probability
one for all $t\geq0$
\begin{equation}
                                        \label{9.7.3}
\tilde x^{(k')}_{t}=x^{(k')}+\int_{0}^{t}
\sqrt{ a^{(k')}(\tilde x^{(k')}_{s})}\,d\tilde w^{(k')}_{s}
+\int_{0}^{t}b^{(k')}(\tilde x^{(k')}_{s}) \,ds.
\end{equation}
Observe  that, in light of Theorem
\ref{theorem 9.7.1}, as for $\tilde x^{(k')}_{t}$,
we have that for any, first, continuous
and, hence, for all Borel nonnegative $f$,
$$
E \int_{0}^{T} 
f(\tilde   x _{t})\,dt\leq N 
T^{ 1/2 }\|\Psi^{-1}_{1/T}f\| _{L_{d}(\bR^{d})}.
$$

Then one passes to the limit in the first term on the right
in \eqref{9.7.3} by literally repeating the 
corresponding part of the proof of Theorem 2.6.1 of \cite{Kr_77}.  

In what concerns the second term, it suffices to observe that
for any $k_{0}$ 
$$
\nlimsup_{k'\to\infty}E\int_{0}^{T}|b^{(k')}(\tilde x^{(k')}_{t})
-b (\tilde  x _{t})|\,dt\leq
\nlimsup_{k'\to\infty}E\int_{0}^{T}|b^{(k')}(\tilde x^{(k')}_{t})
-b^{(k_{0})}(\tilde x^{(k')}_{t})|\,dt
$$
$$
+E\int_{0}^{T}|b^{(k_{0})}(\tilde  x _{t})-b (\tilde  x _{t})|\,dt
\leq N\nlimsup_{n'\to\infty}\|b^{(k')}-b^{(k_{0})}\|_{L_{d}(\bR^{d})}
$$
$$
=N\|b -b^{(k_{0})}\|_{L_{d}(\bR^{d})},
$$
where the constants $N$ are independent of $k_{0}$,
which after sending $k_{0}\to\infty$ shows that the first expression
is zero and this, as in the proof 
of Theorem 2.6.1 of \cite{Kr_77}, allows us to finish
the present proof. The theorem is proved.\qed

Introduce $\cL(\delta,\|\gb\|)$ as the set of
operators
$$
L=(1/2)a^{ij}D_{ij}+b^{i}D_{i},
$$
where $(a^{ij})$ is an $\bS_{\delta}$-valued Borel
function on $\bR^{d}$ and $b=(b^{i})$ is an $\bR^{d}$-valued
Borel
function on $\bR^{d}$ such that $\|b\|_{L_{d}(\bR^{d})}\leq \|\gb\|$.

Here is a generalization for uniformly nondegenerate
operators of the famous  
 Lemma 8 of Aleksandrov \cite{Al_63}
for functions in $W^{2}_{p}$ with $p$ that could be~$<d$.
This result for bounded $b$
is found in \cite{Ca_95} and for
  $b\in L_{d+\varepsilon}(\Omega)$
 in \cite{Fo_98}.
\begin{corollary}
                      \label{corollary 10.7.1}
Let $p\in[d_{0},\infty)$, $D$ be a bounded domain in $\bR^{d}$, and
$u\in W^{2}_{p,\loc}(D)\cap C(\bar D)$.
Let $c$ be a nonnegative measurable function on $D$
and let $L\in \cL(\delta,\|\gb\|)$.
Then in $D$
\begin{equation}
                                                 \label{10.8.1}
u\leq N\|(Lu-cu)_{-}\|_{L_{p}(D)}+
\sup_{\partial D}u_{+},
\end{equation}
where $N$ depends only on $p,d,\delta, \|\gb\|$,
and the diameter of $D$.

\end{corollary}

Proof. It suffices to prove \eqref{10.8.1} only in
$D'=D\cap\{u>0\}$ assuming that this domain is not empty.
Then \eqref{10.8.1} will become stronger if we replace
$D$ with $D'$. Also observe that  on $D'$ we have $(Lu)_{-}\leq
(Lu-cu)_{-}$ since $cu\geq0$. It follows that it suffices
to prove \eqref{10.8.1} in case that $u\geq0$ and $c\equiv0$
on $D$. Having in mind obvious approximation of $D$ from
inside with smooth domains and extending $u$
outside of approximating domains, we may also assume that $D$ is smooth
and $u\in W^{2}_{p }(\bR^{d})$. 

In that case, by Theorem \ref{theorem 9.6.4}, for any $x \in D$
we can find a solution $x_{t}$ of equation \eqref{11.29.2}.
Below in the proof by $x_{t}$ we mean this solution.
In light of Theorem \ref{theorem 9.6.2}, for any $T\in(0,\infty)$,
\begin{equation}
                                                 \label{10.8.3}
u(x )=-E\int_{0}^{\tau\wedge T}Lu(x_{t})\,dt 
+Eu(x_{\tau\wedge T}),
\end{equation}
where $\tau$ is the first exit time of $x_{t}$ from $D$.
Since $E\tau<\infty$ and (cf.~\eqref{9.5.4})
$$
E\int_{0}^{\tau }I_{D}
(x_{t})|Lu(x_{t})|\,dt\leq N\|Lu\|_{L_{d_{0}}(D)}<\infty,
$$
we can pass to the limit in \eqref{10.8.3} as $T\to\infty$
and obtain
$$
u(x )=-E\int_{0}^{\tau }Lu(x_{t})\,dt 
+Eu(x_{\tau })\leq E\int_{0}^{\tau }(Lu(x_{t}))_{-}\,dt
+\sup_{\partial D}u.
$$
After that it only remains to use \eqref{9.5.4} again. \qed

Here is another consequence of It\^o's formula and our previous  
results for elliptic equations. Such results play a crucial role
in the theory of {\em fully nonlinear\/} elliptic equation
providing a tool allowing to pass to the limit
under the sign of a nonlinear operator when there is no convergence
of the derivatives of the functions to which the operator is applied
(see, for instance, Section 4.2 in \cite{Kr_18}).

The following result for $p\geq d$
and $R=\infty$ is obtained in \cite{Kr_12},
however,
with $N$ in \eqref{10.14.4} depending on how fast $\|(|b|-\mu)_{+}\|
_{L_{d}(\bR^{d})}\to 0$ as $\mu\to\infty$.
\begin{theorem}
                                   \label{theorem 10.14.2}
Let $p\geq d_{0}$, $R\in(0,\infty]$, and $L\in\cL(\delta,\|\gb\|)$. Then there exists
a constant
$N=N(p,d,\delta,\|\gb\|)\geq0$ such that for any $\lambda>0$    
and $u\in W^{2}_{p,\loc}(B_{R})\cap C(\bar B_{R})$ 
($B_{\infty}=\bR^{d}$, $C(\bR^{d})$
is the set of bounded continuous functions on $\bR^{d}$)
 we have
\begin{equation}
                                           \label{10.14.4}
\lambda\|u_{+}\|_{L_{p}(B_{R/2})}
\leq N\|(\lambda u-Lu)_{+}\|_{L_{p}(B_{R})}
+N\lambda R^{d/p}e^{-R\sqrt{\lambda}/N}\sup_{\partial B_{R}}u_{+},
\end{equation}
where the last term should be dropped if $R=\infty$.
\end{theorem}

For $p\geq d$ this theorem is proved in \cite{Kr_19}
(see Theorem 3.1 there).
The proof from \cite{Kr_19} carries over to our present situation
almost word for word.

\mysection{Estimates of the time spent in
sets of small measure}

                                             \label{section 10.15.1}

Here is the first result of this section, which will be proved
after some discussion.

\begin{theorem}
                       \label{theorem 9.30.1}

For any 
$\kappa\in(0,1)$ there exist constants $\mu\geq1$ and
$N$, depending only on $ d,\delta,\|\gb\|$, and $\kappa$,
 such that, for any $R\in(0,\infty)$,
$x\in   B_{\kappa R}$, and Borel set
$\Gamma\subset  B_{R}$, the expected time that $x+x_{t}$
spends in $\Gamma$ before exiting from $B_{R}$
is greater than or equal to 
$N^{-1}R^{2}(|\Gamma|/|B_{R}|)^{\mu}$$:$
\begin{equation}
                               \label{10.2.1}
E\int_{0}^{\tau_{R}(x)}I_{\Gamma}(x+x_{t})\,dt
\geq N^{-1}R^{2}(|\Gamma|/|B_{R}|)^{\mu}.
\end{equation}

\end{theorem}

  The results of the kind which follows
are commonly used while establishing the H\"older
continuity of harmonic functions for diffusion processes
or elliptic operators.

\begin{corollary}
                                             \label{corollary 10.11.1}
For any 
$\kappa\in(0,1)$ there exists a constant
$N=N(d,\delta,\|\gb\|,\kappa)$  such that, for any $R\in(0,\infty)$,
$x\in   B_{\kappa R}$, and closed set
$\Gamma\subset  B_{R}$, the probability that $x+x_{t}$
reaches $\Gamma$ before exiting from $B_{R}$
is greater than or equal to 
$N^{-1} (|\Gamma|/|B_{R}|)^{\mu-1/d}$$:$
\begin{equation}
                               \label{10.2.10}
P(\tau_{\Gamma}(x)<\tau_{R}(x))
\geq N^{-1} (|\Gamma|/|B_{R}|)^{\mu-1/d},
\end{equation}
where $\tau_{\Gamma}(x)$ is the first time $x+x_{t}$
hits $\Gamma$
and $\mu$ is taken from Theorem \ref{theorem 9.30.1}.
\end{corollary}

Indeed, set $\gamma=|\Gamma|/|B_{R}|$ and observe
that by Theorem \ref{theorem 9.30.1} and the  conditional
version of Theorem 1.1 of \cite{Kr_19} 
$$
R^{2}\gamma^{\mu}\leq NE
I_{\tau_{\Gamma}<\tau_{R}(x)}E\Big[
\int_{\tau_{\Gamma}}^{\tau_{R}(x)}I_{\Gamma}(x+x_{t})\,dt
\mid \cF_{\tau_{\Gamma}}\Big]
$$
$$
\leq NE
I_{\tau_{\Gamma}<\tau_{R}(x)}R|\Gamma|^{1/d}
=NR^{2}\gamma^{1/d}P(\tau_{\Gamma}(x)<\tau_{R}(x)).
$$

One more consequence of Theorem \ref{theorem 9.30.1}
is the following which will be used in a subsequent article
to show that diffusion processes with drift in $L_{d}$
are strong Markov with strong Feller semigroup
(see \cite{Kr_20}).

\begin{corollary}
                                    \label{corollary 11.1.1}
For $R,\kappa,x$, and $\Gamma$ as in the theorem
set
$$
\phi_{t}(\Gamma)(x)=\int_{0}^{t}I_{\Gamma}(x+x_{s})\,ds.
$$

Then there exists $\theta>0$, depending only on $d,\delta,\|\gb\|$,
and $\kappa$, such that
\begin{equation}
                                          \label{11.5.4}
P (\phi_{\tau_{R}(x)}(\Gamma)(x)\geq\theta\gamma^{\mu} R^{2})
\geq N^{-1}\gamma^{2\mu},
\end{equation}
where $\gamma=  |\Gamma|/|B_{R}|$, $N=N(d,\delta,\|\gb\|,\kappa)$, and $\mu$ is the
same as in Theorem \ref{theorem 9.30.1}.

\end{corollary}

Indeed, \eqref{10.2.1} implies that for any $\theta>0$ and $\phi:=
\phi_{\tau_{R}(x)}(\Gamma)(x)$,
$$
N^{-1}_{1}R^{2}\gamma^{\mu}
\leq \theta\gamma^{\mu} R^{2}P(\phi\leq\theta\gamma^{\mu} R^{2})+
E \phi I_{\phi\geq\theta \gamma^{\mu}R^{2}}
$$
$$
\leq \theta \gamma^{\mu}R^{2}+P^{1/2}(\phi\geq\theta\gamma^{\mu} R^{2})
 E^{1/2}\tau^{2}_{R}(x)
$$
$$
\leq \theta\gamma^{\mu} R^{2}+N_{2}R^{2}P^{1/2}(\phi
\geq\theta \gamma^{\mu} R^{2}),
$$
where the last inequality is obtained by Lemma 
2.6 of \cite{Kr_19}.
We get
\eqref{11.5.4} for $\theta=(1/2)N^{-1}_{1}  $.

\begin{corollary}
                      \label{corollary 10.1.1}
For any $R\in(0,\infty),\kappa\in(0,1)$,
$x\in B_{\kappa R}$, and Borel nonnegative $f$
$$
\int_{B_{R}}f^{1/(2\mu)}(y)\,dy\leq NR^{d-1/\mu}
\Big(E\int_{0}^{\tau_{R}(x)}f(x+x_{t})\,dt\Big)^{1/(2\mu)},
$$
where $N=N(d,\delta,\|\gb\|,\kappa)$.

\end{corollary}

Indeed, without losing generality assuming that
$f=0$ outside $B_{R}$ and
setting
$$
u:=E\int_{0}^{\tau_{R}(x)}f(x+x_{t})\,dt,
$$
we have that 
for any $\lambda>0$
$$
u\geq \lambda
E\int_{0}^{\tau_{R}}I_{f(x+x_{t})\geq\lambda}\,dt
\geq\lambda N^{-1}R^{2}\big(|\{f\geq\lambda\}|/|B_{R}|\big)^{\mu},
$$
$$
|\{f\geq\lambda\}|\leq N R^{-2/\mu}\lambda^{-1/\mu}
|B_{R}|u^{1/\mu}.
$$
It follows that for any $c>0$
$$
\int_{B_{R}}f^{1/(2\mu)}(y)I_{f>c}\,dy=
\big(1/(2\mu)\big)\int_{c}^{\infty}\lambda^{1/(2\mu)-1}
|\{f(y)>\lambda\}|\,d\lambda
$$
$$
\leq N R^{-2/\mu} 
|B_{R}|u^{1/\mu}c^{-1/(2\mu)}.
$$
Also
$$
\int_{B_{R}}f^{1/(2\mu)}(y)I_{f\leq c}\,dy\leq
c^{1/(2\mu)}|B_{R}|.
$$
For $c=uR^{-2}$ we have
$$
R^{-2/\mu} 
|B_{R}|u^{1/\mu}c^{-1/(2\mu)}=
c^{1/(2\mu)}|B_{R}|,
$$
$$
\int_{B_{R}}f^{1/(2\mu)}(y)\,dy
\leq N u^{1/(2\mu)}R^{d-1/\mu}.
$$
This is what is claimed.

Another corollary is a generalization
of the Fanghua Lin estimate for operators
with summable drift which is one of the main tools in
the Sobolev space theory of fully nonlinear elliptic
equations (see, for instance, \cite{Kr_18}
and \cite{Kr_20_1}). 

\begin{theorem}
                       \label{theorem 10.6.10}
Let  $R\in(0,\infty)$,  $p\in[d_{0},\infty)$, 
$u\in W^{2}_{p,\loc}(B_{R})\cap C(\bar B_{R})$,
$L\in \cL(\delta,\|\gb\|)$, and
$c\in L_{d_{0}}(B_{R})$,
$c\geq0$. Then
\begin{equation}
                              \label{10.6.1}
\Big(\dashint_{B_{R}}|D^{2}u|^{1/(2\mu)} 
\,dx\Big)^{2\mu} \leq N \Big(\dashint_{B_{R}}|L u-cu|^{p}\,dx \Big)^{1/p}
 +NR^{ -2}\sup_{\partial B_{R}}|u|,
\end{equation}
where $\mu$ is taken from Theorem \ref{theorem 9.30.1} 
with $\kappa=1/2$
and $N$ depends only on
$d,\delta,\|\gb\|,p$, and $R^{2-d/d_{0}}\|c\|_{L_{d_{0}}(B_{R} )}$.
\end{theorem}

Proof. H\"older's inequality allows us to
concentrate on the case of $p=d_{0}$.
On the account of moving $R$,
we may assume that 
$u\in W^{2}_{p}(B_{R})$ and then, by using scaling, that $R=1$.
After that we observe that 
$$
\|Lu\|_{L_{p}(B_{1})}\leq\|Lu-cu\|_{L_{p}(B_{1})}
+\| c\|_{L_{p}(B_{1})}\sup_{B_{1}}|u|
$$
$$
\leq \big(1+N\| c\|_{L_{p}(B_{1})}\big)\|Lu-cu\|_{L_{p}(B_{1})}
+\| c\|_{L_{p}(B_{1})}\sup_{\partial B_{1}}|u|
$$
and reduce the case of general $c$ to the one with $c\equiv0$.
In that case,
it is easy to see that for sufficiently
small $\varepsilon=\varepsilon(d,\delta)>0$
there is an operator $L'\in\cL(\delta/2,\|\gb\|)$
such that for our function $u$ we have
$$
L'u=Lu+\varepsilon|D^{2}u|.
$$
Then, if $x'_{t}$
is the process corresponding to $L'$ and starting at 
the origin by It\^o's formula we get (cf.~the proof
of Corollary \ref{corollary 10.7.1})
$$
\varepsilon 
E\int_{0}^{\tau}|D^{2}u|(x'_{t})\,dt=-
u(0)-E\int_{0}^{\tau}Lu(x'_{t})\,dt+Eu(x'_{\tau}),
$$
where $\tau$ is the first exit time of $x'_{t}$
from $B_{1}$. After that it only remains to use
Corollary \ref{corollary 10.1.1}, 
Corollary \ref{corollary 10.7.1}, and 
Theorem \ref{theorem 9.5.1}.\qed  

\begin{remark}
It is standard that Corollary \ref{corollary 10.1.1}
 implies not only \eqref{10.6.1} but also a similar estimate
in half balls and similar estimates  
for $|Du|$
(see, for instance, Section 9.4 in \cite{Kr_18}).  
\end{remark}
 
To prove Theorem \ref{theorem 9.30.1} we need 
three lemmas. In their proofs we, actually,
 translate into probability
language the arguments 
from \cite{KS_80}
reproduced, for instance, in \cite{Kr_18}.
In turn, the arguments in \cite{KS_80}
have their origin in
\cite{KS_79} written 
in the probabilistic language. It is worth noting that
some of  the arguments in \cite{KS_80}
are rewritten in the probabilistic language
in \cite{Gr_84} when $b$ is bounded and $x_{t}$
is a {\em solution\/} of a stochastic equation
with nonrandom regular coefficients.
We start with the following. 

\begin{lemma}
                           \label{lemma 9.28.1}
Let $\kappa\in(0,1)$. Then there is a constant $\xi=\xi(\kappa,
d,\delta,\|\gb\|)
\in(0,1)$ such that for any $R\in(0,\infty)$,
 Borel set $\Gamma\subset B_{R}$
satisfying $|\Gamma|\geq \xi|B_{R}|$,
and $x\in B_{\kappa R}$ we have
\begin{equation}
                               \label{9.28.1}
E\int_{0}^{\tau_{R}(x)}I_{\Gamma}(x+x_{t})\,dt
\geq \nu E\tau_{R}(x),
\end{equation}
where   $\nu=\nu(d,\delta,\|\gb\|,
\kappa)\in(0,1)$.
\end{lemma}

Proof. Fix $x$ with $|x|\leq\kappa R$ and define $\gamma$ as the first exit time of
$x+x_{t}$ from $B_{2R}(x)$. By Corollary 2.12 of \cite{Kr_19} and
Theorem \ref{theorem 9.5.1} 
$$
E\tau_{R}(x)-E\int_{0}^{\tau_{R}(x)}I_{\Gamma}(x+x_{t})\,dt=
E\int_{0}^{\tau_{R}(x)}I_{B_{R}\setminus\Gamma}(x+x_{t})\,dt
$$
$$
\leq E\int_{0}^{\gamma}I_{B_{R}\setminus\Gamma}(x+x_{t})\,dt
\leq NR|B_{R}\setminus\Gamma|^{1/d}
$$
$$
=NR^{2}\big(1-|\Gamma|/|B_{R}|\big)^{1/d}
\leq NE\tau_{R}(x)(1-\xi)^{1/d},
$$
where the constants $N$ depend only on $\kappa,d,\delta$, and $\|\gb\|$. We see how to choose
$\xi$ to satisfy \eqref{9.28.1}
with a $\nu=\nu(d,\delta,\|\gb\|,
\kappa)\in(0,1)$.
The lemma is proved. \qed

Next, we need a fact from the geometric measure theory.
Take $R\in(0,\infty)$, $\zeta \in(0,1)$, and a Borel set
$\Gamma\subset B_{R}$ such that $|\Gamma|<\zeta|B_{R}|$.
Define $\frA$ as the collection
of all open balls $B\subset B_{R}$ such that
\begin{equation}
                                                \label{9.16.1}
|\Gamma\cap B|\geq \zeta  |B|.
\end{equation}

Observe that this collection is nonempty   
and the union of its elements 
$$
\Gamma'=\bigcup_{B \in\frA}B.
$$
contains almost any
point of $\Gamma$ since
almost  any
point of $\Gamma$ is its density point.
Also 
observe that $\Gamma'$ is an open set, since all the $B$'s
are open.

Then for $\varepsilon\in(0,1)$
denote by $\frA_{\varepsilon}$ the set of
$B\in \frA$ such that $|B|\geq\varepsilon$.
Finally, recall that if $B$ is an open ball and $\kappa\in(0,1)$,
we write $\kappa B$ for the concentric open ball of
radius $\kappa$ times that of $B$
and set
$$
\Gamma'_{\kappa }=\bigcup_{B\in\frA }
\kappa B,\quad
\Gamma'_{\kappa,\varepsilon}=\bigcup_{B\in\frA_{\varepsilon}}
\kappa B.
$$

\begin{lemma}
                                             \label{lemma 9.16.1}
1. We have $|\Gamma\setminus\Gamma'|=0$
and
$$
|\Gamma'|\geq\Big(1+\frac{1-\zeta}{3^{d}}\Big)
|\Gamma|.
$$

2. There exists $\kappa=\kappa(d,\zeta)\in(0,1)$
and $\theta=\theta(d,\zeta)>1$ such that
for all sufficiently small $\varepsilon>0$
there exists a  closed set $\Gamma''_{\varepsilon}
\subset \Gamma'_{\kappa,\varepsilon}$ such that
$$
|\Gamma''_{\varepsilon}|\geq\theta |\Gamma|.
$$

\end{lemma}

Proof. The first assertion, a parabolic version of which is found in
\cite{KS_80}, is 
proved in Lemma 1.1 of \cite{Sa_80}.
To prove the second one it suffices to
observe that, obviously, $\Gamma'_{\kappa,\varepsilon}
\uparrow \Gamma'_{\kappa}$ as
$\varepsilon\downarrow0$ and, 
similarly to Lemma 2.4 of \cite{KS_80},
$| \Gamma'_{\kappa}|\ge \kappa^{d}|\Gamma'|$. 
The lemma is proved. \qed

Since $\Gamma''_{\varepsilon}
\subset \Gamma'_{\kappa,\varepsilon}$,
the closed set $\Gamma''_{\varepsilon}$ is covered by the family $\{\kappa B:
B\in\frA_{\varepsilon}\}$. Then there is
finitely many $B(1),...,B(n)
\in\frA_{\varepsilon}$ such that
$$
\Gamma''_{\varepsilon}
\subset \bigcup_{i=1}^{n}\kappa B(i)
=:\Pi_{\varepsilon}.
$$
Next, for $x\in \Pi_{\varepsilon}$ define
$i(x)$
as the first $i\in\{1,...,n\}$ for which
$x\in \kappa B(i)$.  Also set $B(0)=B_{R}$
and $i(x)=0$ if $x\in\partial B_{R}$.
Now define recursively $\gamma^{0}=0$,
$\tau^{1}$ as the first time after $\gamma^{0}$ when $x_{t}$ exits
from $B_{R}\setminus \Gamma''_{\varepsilon}$,
$\gamma^{1}$ as the first  time after $\tau^{1}$
when $x_{t}$ exits from $B(i(x_{\tau^{1}}))$,
and generally, for $n=2,3,...$ define
$\tau^{n}$ as the first time after $\gamma^{n-1}$ when $x_{t}$ exits
from $B_{R}\setminus \Gamma''_{\varepsilon}$,
$\gamma^{n}$ as the first  time after $\tau^{n}$
when $x_{t}$ exits from $B(i(x_{\tau^{n}}))$.
It is easy to check that so defined
$\tau^{n}$ and $\gamma^{n}$ are stopping times
and, since $|B(i)|\geq\varepsilon$ and the trajectories of $x_{t}$ are continuous,
$\tau^{n}\uparrow \tau_{R}$ as $n\to\infty$.
Furthermore, since (a.s.) $\tau_{R}$ is finite, (a.s.) all the $\tau^{n}$'s equal $\tau_{R}$ for all large $n$.

Now in the above general constructions we set
$\zeta=\xi$, where $\xi$ is taken from
Lemma \ref{lemma 9.28.1}.
\begin{lemma}
                           \label{lemma 9.28.2}
Set $\zeta=\xi$, where $\xi$ is taken from
Lemma \ref{lemma 9.28.1}. Then for $x$
such that $|x|\leq\kappa R$ we have
\begin{equation}
                            \label{9.28.3}
E\int_{0}^{\tau_{R}(x)}I_{\Gamma}(x+x_{t})\,dt
\geq \nu  
E\int_{0}^{\tau_{R}(x)}I_{\Gamma''_{\varepsilon}}(x+x_{t})\,dt,
\end{equation}
where $\nu$
is taken from
Lemma \ref{lemma 9.28.1}.
\end{lemma}

Proof. By the conditional  version of 
Lemma \ref{lemma 9.28.1} (a.s.)
$$
E\Big[\int_{\tau^{k}}^{\gamma^{k}}I_{\Gamma''_{\varepsilon}}(x+x_{t})\,dt
\mid \cF_{\tau^{k}}\Big]
\leq
E\Big[\int_{\tau^{k}}^{\gamma^{k}}I_{B(i(x_{\tau^{k}}))}(x+x_{t})\,dt
\mid \cF_{\tau^{k}}\Big]
$$
$$
\leq \nu^{-1}
E\Big[\int_{\tau^{k}}^{\gamma^{k}}I_{\Gamma\cap B(i(x_{\tau^{k}}))}(x+x_{t})\,dt
\mid \cF_{\tau^{k}}\Big]
\leq \nu^{-1}
E\Big[\int_{\tau^{k}}^{\gamma^{k}}I_{\Gamma}(x+x_{t})\,dt
\mid \cF_{\tau^{k}}\Big].
$$

Hence,
$$
E\int_{0}^{\tau_{R}(x)}I_{\Gamma''_{\varepsilon}}(x+x_{t})\,dt=\sum_{k=1}^{\infty}
E\int_{\tau^{k}}^{\gamma^{k}}I_{\Gamma''_{\varepsilon}}(x+x_{t})\,dt
$$
$$
\leq\nu^{-1}\sum_{k=1}^{\infty}
E\int_{\tau^{k}}^{\gamma^{k}}I_{\Gamma}(x+x_{t})\,dt 
\leq\nu^{-1} E\int_{0}^{\tau_{R}(x)}I_{\Gamma}(x+x_{t})\,dt.
$$
This proves the lemma.\qed

{\bf Proof of Theorem \ref{theorem 9.30.1}}.
Take $\xi$ from Lemma \ref{lemma 9.28.1}.
The results of Lemma's \ref{lemma 9.16.1} and \ref{lemma 9.28.2} can
 be summarized as follows: For $x\in \kappa B_{R}$, if a measurable $\Gamma\subset B_{R}$
is such that $|\Gamma|<\xi|B_{R}|$,
then there exists a closed set
$\Gamma_{1}\subset B_{R}$ such that
$|\Gamma_{1}|\geq \theta|\Gamma|$, where
$\theta=\theta(d,\xi)>1$ is taken from Lemma
\ref{lemma 9.16.1}, and
$$
E\int_{0}^{\tau_{R}(x)}I_{\Gamma}(x+x_{t})\,dt
\geq \nu  
E\int_{0}^{\tau_{R}(x)}I_{\Gamma_{1}}(x+x_{t})\,dt
$$
with $\nu$ from Lemma \ref{lemma 9.28.2}.
Naturally, if $|\Gamma_{1}|<\xi|B_{R}|$, which
only happens if $|\Gamma|\leq (\xi/\theta)|B_{R}|$,
then there exists a closed set
$\Gamma_{2}\subset B_{R}$ such that
$|\Gamma_{2}|\geq \theta|\Gamma_{1}|\geq \theta^{2}
|\Gamma|$ and
$$
E\int_{0}^{\tau_{R}(x)}I_{\Gamma}(x+x_{t})\,dt
\geq \nu 
E\int_{0}^{\tau_{R}(x)}I_{\Gamma_{1}}(x+x_{t})\,dt
\geq \nu^{2}  
E\int_{0}^{\tau_{R}(x)}I_{\Gamma_{2}}(x+x_{t})\,dt.
$$
We continue in a natural way and see that,
if $n$ is such that $|\Gamma_{n}|<\xi|B_{R}|$, which
only happens if $|\Gamma|\leq (\xi/\theta^{n})|B_{R}|$,
then there exists a closed set
$\Gamma_{n+1}\subset B_{R}$ such that
$|\Gamma_{n+1}|\geq \theta^{n+1}
|\Gamma|$ and
$$
E\int_{0}^{\tau_{R}(x)}I_{\Gamma}(x+x_{t})\,dt
\geq \nu^{n+1}  
E\int_{0}^{\tau_{R}(x)}I_{\Gamma_{n+1}}(x+x_{t})\,dt.
$$
Let $n_{0}$ be the largest $n$ for which the construction
of $\Gamma_{n+1}$ is still possible, that is $|\Gamma_{n_{0}}|<\xi|B_{R}|$ and
$|\Gamma_{n_{0}+1}|\geq\xi|B_{R}|$.
Since $|\Gamma_{n}|\geq \theta^{n}
|\Gamma|$, we have 
$$
n_{0}\leq\Big\lfloor\big(\ln\big(|B_{R}|/|\Gamma|\big)/\ln\theta\Big\rfloor.
$$
Since by Lemma \ref{lemma 9.28.1}
$$
E\int_{0}^{\tau_{R}(x)}I_{\Gamma_{n_{0}+1}}(x+x_{t})
\,dt\geq \nu E\tau_{R}(x),
$$
we have
$$
E\int_{0}^{\tau_{R}(x)}I_{\Gamma}(x+x_{t})\,dt
\geq \nu^{n_{0}+2} E\tau_{R}(x).
$$
We take into account that by Corollary 2.12 of \cite{Kr_19}
  $E\tau_{R}(x)\geq N_{1}^{-1}R^{2}$
and come to \eqref{10.2.1} with   
$$
\mu=-(\ln \nu)/\ln\theta,\quad N=N_{1}\nu^{-3}.
$$

This takes care of the case in which
$|\Gamma|<\xi|B_{R}|$. To include the case
$|\Gamma|\geq\xi|B_{R}|$, in light of Lemma
\ref{lemma 9.28.1}, it suffices to increase
the above $N$ in an obvious way. The theorem
is proved. \qed

\begin{corollary}
                                       \label{corollary 10.21.1}
 
Let $R\in(0,\infty)$, $\gamma\in(0,1)$,
 and assume that a closed set $\Gamma\subset
B_{R}$ is such that, for any $r\in(0,R)$, $|B_{r}\cap\Gamma|\geq
\gamma |B_{r}|$. Then there exist    constants $\alpha\in(0,1)$
and $N$, depending only on $d,\delta,
\|\gb\|$, and $\gamma$, such that, for any $ x\in B_{R/2}$,
\begin{equation}
                                     \label{10.21.1}
 P(\tau_{R}(x)<\tau_{\Gamma}(x))\leq N(x/R)^{\alpha}.
\end{equation}
\end{corollary}

Indeed, let $R^{n}=R2^{-n}$, $\Gamma_{n}=\Gamma\cap B_{R_{n}}$,
and $A_{n}=\{\tau_{R_{n}}(x)<\tau_{\Gamma_{n}}(x)\}$, $n=0,1,...$.
Then
by Corollary \ref{corollary 10.11.1}, for $|x|\leq R_{n+1}$,
 $
P(A_{n})\leq q = q(d,\delta,\|\gb\|,\gamma)<1
 $.
The conditional  version of this says that on the set
$A_{n+1}$ we have (a.s.)
 $
P(A_{n}\mid \cF_{\tau_{R_{n+1}}})\leq q
 $.
Since, for $|x|\leq R/2$, $A_{0}=\bigcap_{n=0}^{n(x)}A_{n}$, where
$n(x)= \lfloor\ln_{2}(R/x)\rfloor-1 $ ($\geq 0$), we have
$$
P(A_{0})\leq q^{n(x)+1}\leq q^{-1}(x/R)^{\ln_{2}(1/q)},
$$
which is just a different form of \eqref{10.21.1}.

The following result will be used in a subsequent article
on fully nonlinear elliptic equations with singular lower 
order terms (see \cite{Kr_20_1}).

 \begin{theorem}
                                        \label{theorem 10.20.3}
Let $D$ be a bounded domain in $\bR^{d}$, $0\in\partial D$,
and assume that for some constants $\rho,\gamma>0$ and 
any $r\in (0,\rho)$ we have $|B_{r}\cap D^{c}|\geq \gamma
|B_{r}|$. Then there exists $\beta=\beta(d,\delta,\|\gb\|,\gamma  )>0$
such that, for any nonnegative $f\in L_{d_{0}}(D)$ and $x\in D$,
\begin{equation}
                                                 \label{10.20.4}
u(x):=E\int_{0}^{\tau(x)}f(x+x_{t})\,dt
\leq N|x|^{\beta}\|f\|_{L_{d_{0}}(D)},
\end{equation}
where $\tau(x)$ is the first exit time of $x+x_{t}$
from $D$ and $N$ depends only on $ d,\delta,\|\gb\|,\gamma,\rho$, and the diameter
of $D$.
\end{theorem}

Proof. In light of Theorem \ref{theorem 9.5.1}
we may concentrate on $x\in B_{\rho/2}$ with $|x|\leq1$. The 
conditional 
version of this theorem allows us to write that,
for $2|x|\leq r<\rho$  and $\tau^{r}(x)$ being the first
exit time of $x+x_{t}$ from $B_{r}\cap D$,
$$
u(x)=E\int_{0}^{\tau^{r}(x)}f(x+x_{t})\,dt
+EI_{\tau^{r}(x)<\tau(x)}E\Big(\int_{\tau^{r}(x)}^{\tau(x)}f(x+x_{t})\,dt\mid
\cF_{\tau_{r}(x)}\Big)
$$
$$
\leq N r^{2-d/(2p)}\|f\|_{L_{d_{0}}(D)}
+N\|f\|_{L_{d_{0}}(D)}P(\tau^{r}(x)<\tau(x)).
$$
Observe that $\{\tau^{r}(x)<\tau(x)\}\subset\{\tau^{r}(x)<\tau_{\Gamma_{r}}(x)\}$,
where $\Gamma_{r}=\bar B_{r}\cap D^{c}$, and by Corollary
\ref{corollary 10.21.1} we have
$P(\tau^{r}(x)<\tau(x))\leq N(x/r)^{\alpha}$. Thus,
for any $r\in [2|x|,\rho)$
$$
u(x)\leq N\|f\|_{L_{p}(D)}\big( r^{2-d/(2d_{0})}
+ (x/r)^{\alpha}\big).
$$
 By choosing $r=\sqrt{2|x|\rho}$,
we get the result with $\beta=\alpha/2$
since $1-d/(4d_{0})>1/2>\alpha/2$. The theorem is proved.
\qed
\begin{remark}
                                       \label{remark 10.20.3}
If $b$ and $f$ are bounded and a part of $\partial D$ near the origin
is flat, then one can take $\beta=1$ in \eqref{10.20.4}.
However, even in the case of flat boundary and bounded $f$,
if $b\in L_{d}$, then in the general case 
certainly $\beta<1$ (see Example 4.1  in \cite{Sa_10}) and most likely $\beta\to0$
as $\delta\to0$.

\end{remark}
\begin{corollary}
                                       \label{corollary 1.6.1}
Under the assumptions of Theorem \ref{theorem 10.20.3}
suppose that we are given a function $u\in W^{2}_{d_{0},\loc}
(D)\cap  C(\bar D)$. 
Let $w(r)$ be a concave continuous function on $[0,\infty)$
such that $w(0)=0$ and $|u(x)-u(0)|\leq
w(|x|)$ for all $x\in \partial D$.
Then for $x\in D$ we have
\begin{equation}
                                               \label{1.6.4}
|u(x)-u(0)|\leq N|x|^{\beta}\|Lu\|_{L_{d_{0}}(D)}
+\omega\big(N|x|^{\beta/2}),
\end{equation}
where $\beta$ is the same as in Theorem \ref{theorem 10.20.3}
and $N$   depend on the data in the same way as
in Theorem \ref{theorem 10.20.3}.
\end{corollary}

Indeed, define $f=-Lu$. Clearly, we may assume that
$f\in L_{d_{0}}(D)$. Then take a sequence of domains
$D_{n}\subset \bar D_{n}\subset D$ such that
$D_{n}\uparrow D$, denote by $\tau_{n}(x)$
the first exit time of $x+x_{t}$ from $D_{n}$,
and use It\^o's formula (cf.~the proof
of Corollary \ref{corollary 10.7.1})
to conclude that
$$
u(x)=E\int_{0}^{\tau_{n}(x)}f(x+x_{t})\,dt
+Eu(x+x_{\tau_{n}(x)}).
$$
In light of \eqref{9.5.4} we can pass to the limit
as $n\to\infty$ and owing to \eqref{10.20.4} to conclude
$$
u(x)-u(0)=E\int_{0}^{\tau (x)}f(x+x_{t})\,dt
+E[u(x+x_{\tau (x)})-u(0)]
$$
$$
\leq N|x|^{\beta}\|f\|_{L_{d_{0}}(D)}
+Ew(|x+x_{\tau (x)}|).
$$
Here, by Jensen's inequality, the last  term 
is less than $w$ evaluated at the square root
of
$$
E|x+x_{\tau (x)}|^{2}=|x|^{2}+ 2E\int_{0}^{\tau (x)}
[ \tr a_{s}+(x+x_{s})b_{s}]\,ds
$$
$$
\leq
|x|^{2}+N|x|^{\beta}(1+\text{\rm diam}\,
(D)\|\gb\|).
$$
This yields an estimate for $u(x)-u(0)$ from above.
Similarly one estimates it from below.

Next, we study the probability to pass
through narrow tubes, which was, as is mentioned before,
the starting point of the whole theory
presented here.
 We represent the points in $\bR^{d}$ as $x=(x^{1},x')$,
where $x^{1}\in\bR$ and  $x'\in\bR^{d-1}$.

\begin{theorem}
                                         \label{theorem 9.2.1}
Let $\kappa\in[1/2,1)$. Then there exist  $T_{1}>1>T_{0}>0,p_{0}>0$, 
depending only on $\kappa$, $d, \|\gb\|$, and $\delta$,  such that, if  
$R\in(0,\infty)$, $n\in\{2,3,...\}$,
and an open round cylinder $C$ in $\bR^{d}$ with base
being a ball in $\bR^{d-1}$ of radius $R $ and length $nR$
is defined by
$C=(0,nR)\times\{x':|x'|<R\} $, then, for any $x =
(R,x' )$ with $|x' |\leq \kappa R$,
the probability
that $x +x_{t}$ will first exit from the cylinder through
$\{nR\}\times \{x':|x'|<(1-\kappa)R\}$
and this will happen in the time interval $[(n-1)T_{0}R^{2},
(n-1)T_{1}R^{2}]$
is greater than or equal to
$p_{0}^{n-1}$.

\end{theorem} 

 Proof. As usual we may concentrate on $R=1$.
Let us call the sections of $\bar C$
by hyperplanes $x^{1}=k$ disks. We contract them
to their centers with the coefficient of contraction, say $c$
 and call
the results $c$-subdisks. In this terminology we need to estimate
 from below the probability of the event $A$
that our process starting from a point on  the
$\kappa$-subdisk lying at  the distance $1$ from the base
 on which $x^{1}=0$
will first exit from $C$ during the time interval $[nT_{0} ,
nT_{1} ]$ through the (smaller) $(1-\kappa)$-subdisk on the other base.

Let $C_{k}=C\cap\{k < x^{1}<k+2\}$, $k=0,1,...,n-2$.
Then for $A$ to happen it suffices for the process to 
consecutively   exit
from each  
  $C_{k }$, $k=0,1,...,n-2$, 
through the $(1-\kappa)$-subdisk where $x^{1}=k+2$
 during the time interval $[ (k+1)T_{0} ,
(k+1)  T_{1} ]$.
 By using conditional
expectations we easily see that $P(A)\geq p_{0}^{n-1}$,
where $p_{0}$ is the estimate from below in terms of only
$d, \|\gb\|$, and $\delta$ of the probability of the event $A_{0}$
that our process will first exit from $C_{0}$ through the 
$(1-\kappa)$-subdisk
on which $x^{1}=2$ during the time interval $[ T_{0} ,
 T_{1} ]$.

Let $B$ be the open unit ball  centered at $x_{0}=(2-\kappa,0)$,
which is slightly off the center of $C_{0}$. 
Since $|B\cap\{x^{1}\geq 2\}|=N(\kappa)
>0$ and $|(1,x')-x_{0}|^{2} \leq (1-\kappa)^{2}
+\kappa^{2}<1$ if $|x'|\leq \kappa$,
 by Corollary
\ref{corollary 10.11.1} we have $P(A'_{0})\geq 3p_{0}
=3p_{0}(\kappa,d,\delta,\|\gb\|)
>0$, where $A'_{0}$ is the event that $x +x_{t}$
reaches $\overline{B\cap\{x^{1}\geq 2\}}$ before exiting from $B$,
that is reaches $\{x^{1}=2\}\cap\{x':|x'|\leq 1-\kappa\}$
before exiting from $B$. We also know that $E\tau_{B}\leq N$,
where $\tau_{B}$ is the first exit time
of $x +x_{t}$ from $B$. 
Hence, there is $T_{1}>1$ such that
$P(\tau_{B}>T_{1})\leq  p_{0}$.
We also know (see Theorem 2.10 in \cite{Kr_19})
that there exists $T_{0}\in(0, 1)$ such that 
$P(\tau_{B}\leq T_{0})\leq  p_{0}$. Hence,
$P(A'_{0},\tau_{B}\in [T_{0},T_{1}])\geq p_{0}$.
Since $B\cap\{x^{1}<2\}\subset C_{0}$,
 obviously, $A'_{0}\cap\{\tau_{B}\in [T_{0},T_{1}]\}\subset A_{0}$ 
and the theorem is proved.
\qed
 
\begin{corollary}
                                          \label{corollary 10.28.1}
Let $R\in(0,\infty)$ and $|x|\leq R$. Then there is 
a constant $N=N(d,\delta,\|\gb\|)$ such that the expected time spent
by $x_{t}$ in $B_{R}(x)$ before exiting from
$B_{2R}$ is greater than $N^{-1}R^{2}$$:$
$$
R^{2}\leq NE\int_{0}^{\tau_{2R}}I_{B_{R}(x)}(x_{t})\,dt.
$$
 
\end{corollary}

 Indeed, as always we may assume that $R=1$ and then
by Theorem \ref{theorem 9.2.1}
with  probability $p=p(d,\delta,\|\gb\|)>0$ the process $x_{t}$
reaches $\bar B_{1/2}(x)$ before exiting from $B_{2}$.
After that happens the expected time spent in $B_{1}(x)$
before $\tau_{2 }$
is greater than the expected time spent in $B_{1 }(x )$
before exiting from it. Then it only remains 
to use Corollary 2.12 of \cite{Kr_19},
 according to which the expected exit time
from $B_{1}(x )$ starting from a point
in $\bar B_{1/2}(x)$ is greater than $N^{-1} $.

This corollary easily implies the so-called
doubling property of the Green's measure of $x_{t}$.
Let $D\subset \bR^{d}$ be a  bounded domain
containing the origin.
Then the Green's measure of $x_{t}$ in $D$ is defined  
by
$$
G(\Gamma)=E\int_{0}^{\tau}I_{\Gamma}(x_{t})\,dt,
$$
where $\tau$ is the first exit time of $x_{t}$ from $D$.
As we know from the above, $G$ has a density summable to the power
of $d_{0}/(d_{0}-1)$.

\begin{theorem}[doubling property]
                                          \label{theorem 10.28.2}
Let a
ball $B\subset D$ be such that $2B\subset D$. Then
$G(B)\leq NG((1/2)B)$, where $N=N(d,\delta, \|\gb\|)$.
 \end{theorem}

Proof. We may assume that $D$ is connected and the radius
of $B$ is one. Then define  $\tau^{D}$ as the first exit
time of $x_{t}$ from $D$  
and introduce recurrently, for $n=0,1,2,...$, $\gamma^{0}=0$,
$$
\tau^{n }=\inf\{t\geq \gamma^{n}:x_{t}\in \bar B\},\quad
\gamma^{n+1}=\inf\{t\geq\tau^{n}:x_{t}\not\in  (3/2) B\}.
$$
By the conditional versions of, first, Corollary 2.1 of \cite{Kr_19}
and then Corollary \ref{corollary 10.28.1}
we have
$$
G(B)=\sum_{n=1}^{\infty}EI_{\tau^{n}<\tau^{D}}
E\Big(\int_{\tau^{n}}^{\gamma^{n+1}}I_{B}(x_{t})\,dt
\mid \cF_{\tau^{n}}\Big)
$$
$$
\leq N \sum_{n=1}^{\infty}EI_{\tau^{n}<\tau^{D}}
\leq N \sum_{n=1}^{\infty}EI_{\tau^{n}<\tau^{D}}
E\Big(\int_{\tau^{n}}^{\gamma^{n+1}}I_{(1/2)B}(x_{t})\,dt
\mid \cF_{\tau^{n}}\Big)
$$
\begin{equation}
                                                 \label{10.28.7}
=NG((1/2)B),
\end{equation}
which proves the theorem. \qed

\begin{corollary}[$A_{\infty}$-property of $G$]
                                  \label{corollary 10.28.100}
There are constants $\mu\geq1$ and $N$,
depending only on $d,\delta,\|\gb\|$, such that for any
ball $B$ satisfying $2B\subset D$ and Borel $\Gamma\subset B$
we have
\begin{equation}
                                          \label{10.28.5}
N\frac{G(\Gamma)}{G(B)}\geq\Big(\frac{|\Gamma|}{|B|}\Big
)^{\mu}.
\end{equation}
\end{corollary}

Proof. Take the same $\gamma^{n},\tau^{n}$
as in the proof of Theorem \ref{theorem 10.28.2}
and observe that by the conditional  version of
Theorem \ref{theorem 9.30.1} on the set $\{\tau^{n}<\tau^{D}\}$ (a.s.)
$$
E\Big(\int_{\tau^{n}}^{\gamma^{n+1}}I_{\Gamma}(x_{t})\,dt
\mid \cF_{\tau^{n}}
\Big)
\geq N^{-1}R^{2}(|\Gamma|/|B |)^{\mu},
$$
where $R$ is the radius of $B$.
Furthermore,
$$
NR^{2}\geq E\Big(\int_{\tau^{n}}^{\gamma^{n+1}}I_{B}(x_{t})\,dt
\mid \cF_{\tau^{n}}
\Big).
$$
After that it only remains to mimic \eqref{10.28.7}. \qed

Corollary \ref{corollary 10.28.100} is almost identical 
to Corollary 2.3 in \cite{FS_84}. However, there are no lower order terms
in \cite{FS_84} and the comparable situations would be only when
$x_{t}$ were a solution of \eqref{11.29.2}. 
We refer the reader to the proof of Corollary 2.3 in \cite{FS_84}
concerning $A_{\infty}$-weights only pointing out that
Corollary \ref{corollary 10.28.100} is not sufficient for
proving even Theorem \ref{theorem 9.30.1} because
not arbitrary subsets of $B$ could be considered.
On the other hand, $N$ and $\mu$ in \eqref{10.28.5} are independent
of how close $\partial D$ to the origin is
in contrast with Theorem \ref{theorem 9.30.1},
where the starting point of the process is at a distance
at least $(1-\kappa)R$ from the boundary.

\begin{remark}
                                       \label{remark 12.4.1}
Once we know that $G$ is an $A_{\infty}$-weight, it is also
an $A_{p}$-weight for certain large $p$. In particular,
on any closed $\Gamma\subset D$,
$G^{-\alpha}$ is  summable  for some $\alpha>0$.

\end{remark}

Theorem \ref{theorem 9.2.1} allows us to prove
a few more   properties of $x_{t}$. By Corollary 2.7
of \cite{Kr_19} there are  constants $N,\nu>0$, depending only on
$d,\delta$, and $\|\gb\|$, such that 
for any $R, T >0$,
$$
P( \tau_{R}  \geq T)\leq Ne^{-\nu T/R^{2}}.
$$
This turns out to be very close to an optimal result.
\begin{lemma}
                                             \label{lemma 12.8.1}
There are  constants $N,\nu>0$, depending only on
$d,\delta$, and $\|\gb\|$, such that 
for any $R,T>0$,
\begin{equation}
                                                 \label{12.8.1}
NP( \tau_{R} > T 
)\geq e^{-\nu T/R^{2}}.
\end{equation}

\end{lemma}

Proof. We may assume that $R=3$. Then
the cylinder $C=(-1,2)\times \{x':|x'|< 1\}\subset
B_{3}$ and $\tau_{3}> \tau$, where $\tau$
is the first exit
time of $x_{t}$ from $C$. 
Introduce, times of meandering: $\tau_{0}=0$ and, for $n=0,1,2,...$,
let
$\tau_{2n+1}$ be the first exit time of
$x_{t}$ from $(-1,1)\times \{x':|x'|\leq 1\}$ after $\tau_{2n }$,
$\tau_{2n+2}$ be the first exit time of $x_{t}$ 
after $\tau_{2n+1}$   from 
$(0,2)\times \{x':|x'|\leq 1\}$.
Also let $\kappa=1/2$, take $T_{0},T_{1}$ from Theorem
\ref{theorem 9.2.1},   introduce 
$$
A_{2n+1}=\{x_{\tau_{2n+1}}\in \{1\}\times \{x':|x'|\leq 1/2\},
\tau_{2n+1}-\tau_{2n}\in [ T_{0},T_{1}]\},
$$
$$
A_{2n+2}=\{x_{\tau_{2n+2}}\in \{0\}\times \{x':|x'|\leq 1/2\},
\tau_{2n+2}-\tau_{2n+1}\in [ T_{0},T_{1}]\},
$$
and define $n_{0}$ as the least integer   such that 
$n_{0}T_{0}\geq T$. Observe that on the set
$$
\bigcap_{i=1}^{n_{0}}A_{i}
$$
we have
$$
T\leq  n_{0} T_{0}\leq \tau_{n_{0}}<\tau
$$
and hence
$$
P(\tau_{3}>T)\geq P(\tau>T)\geq
P\Big(\bigcap_{k=1}^{n_{0}}A_{k}\big)\geq p_{0}^{n_{0}},
$$
where the last inequality follows from
the conditional version of Theorem
\ref{theorem 9.2.1}. This obviously proves the lemma.  \qed

The following result will be used in a subsequent paper
for establishing Harnack's inequality
for caloric functions related to diffusion processes
with drift in $L_{d}$ (see \cite{Kr_20}).

 \begin{theorem}
             \label{theorem 12.7.2}
Let $R\in(0,\infty)$, $ \kappa,\eta\in(0,1)$,
$x,y\in B_{\kappa R}$, 
 and $\eta^{-1} R^{2}\geq t\geq \eta R^{2}$.
Then there exist $N,\nu>0$, depending only on $ \kappa,\eta,d,\delta$,
and $\|\gb\|$, such that, for any
$\rho\in(0,1]$,
\begin{equation}
                                                 \label{12.9.1}
NP(x+x_{t}\in B_{\rho R}(y),\tau_{   R}(x)> t )\geq \rho^{\nu }.
\end{equation}
\end{theorem}

We prove this theorem after appropriate preparations.

\begin{lemma} 
                    \label{lemma 12.15.1}
If $ \rho_{0}\in(0,1)$, $\xi\in(0,\infty)$, and $\kappa\in[1/2,1)$,
then there exists
$\mu=\mu(d,\delta,\|\gb\|,\kappa, \rho_{0},\xi)>0$ such that
\begin{equation}
                                                 \label{12.7.3}
 P(x+x_{\xi R^{2}}\in B_{\rho_{0}\kappa R }(y),
\tau_{ R}(x)> \xi R^{2}  )\geq \mu,
\end{equation}
whenever $R\in(0,\infty)$,  
$x,y\in B_{\kappa R}$.
\end{lemma}

Proof. While proving \eqref{12.7.3} we may assume that $R=1$.
Then observe that \eqref{12.7.3} becomes stronger if
$\rho_{0}$ becomes smaller. Therefore we may assume that
\begin{equation}
                               \label{12.7.4}
  \rho_{0}\leq \min\big(\kappa^{-1}-1, \xi / T_{1} \big),
\end{equation}

Then also assume, as the first case,
 that $2|x- y|\geq \kappa\rho_{0}$ and 
connect $x$ and $y$ by a round cylinder of length $nr$, where
$$
n=\Big\lfloor\frac{9|x-y|}{\kappa\rho_{0}}\Big\rfloor+1,\quad
r=\frac{|y-x|}{n-1}.
$$
More precisely our cylinder is given by
$$
C=\{x+t(y-x)/|y-x|+ re :t\in(-r,(n-1)r),
e\in \bR^{d},|e|<1,e\perp x\}.
$$
It is not hard to check that, owing to $\rho_{0}\leq\kappa^{-1}-1$,
we have  $C\subset B_{1}$. Also as is easy to see
$$
\kappa\rho_{0}/7\geq r\geq \kappa\rho_{0}/9.
$$

Define
$$
\tau=\inf\{s:x+x_{s}\in \bar B_{ r}(y)  \}.
$$ 
By Theorem \ref{theorem 9.2.1}  we obtain that
with probability not less than $p_{0}^{n-1}$
we have $\tau\leq (n-1)T_{1}r^{2}$ and $\tau<\tau_{1}(x)$.
 Furthermore,
$$
(n-1)T_{1}r^{2}=|y-x|T_{1}r\leq 2T_{1}r\leq T_{1}\rho_{0}\leq\xi.
$$
By Lemma \ref{lemma 12.8.1},
given that $\tau\leq \xi\wedge \tau_{1}(x) $, the probability
that the process $x_{t}$ does not exit
from $B_{ r   }(x_{\tau})$ before time $\xi$,
assuring  that $x+x_{\xi}\in B_{\kappa\rho _{0}}(y)$ and $\tau_{1 }(x)> 
\xi$,
is bigger than
$
N^{-1}e^{-\nu \xi/\rho_{0}^{2}}$. Hence, 
$$
P(x+x_{\eta}\in B_{\kappa\rho_{0} }(y),
\tau_{1 }(x)>\xi)\geq p_{0}^{n-1}N^{-1}
e^{-\nu\xi  /\rho_{0}^{2}}
 \geq N^{-1}e^{-\nu  /\rho_{0}^{2} }=:\mu,
$$
where the last $\nu$ is perhaps  different from the previous one.
This proves \eqref{12.7.3} if $2|x- y|\geq \kappa\rho_{0}$.
If $2|x- y|< \kappa\rho_{0}$ one does not need the first part of the proof.
The lemma is proved. \qed

\begin{lemma} 
                    \label{lemma 12.15.10}
Let $\kappa,\eta\in(0,1)$. Then there are  constants $N,\nu>0$, depending only on
$\kappa,\eta,d,\delta$, and $\|\gb\|$, such that, 
for any $R\in(0,\infty),\rho\in (0,1),  $ and $x \in B_{\kappa R}$,
\begin{equation}
                            \label{12.14.2}
NP\big( \tau_{R}(x) > \eta R^{2},
x+x_{\eta R^{2}}\in B_{\rho R} 
\big) \geq \rho^{\nu }.
\end{equation}
\end{lemma}

Proof. 
We may assume that $\kappa\in[1/2,1)$. Estimate \eqref{12.7.3},
where we take $\xi=\eta$, $y=0$, and $\rho_{0}$ equal to the right-hand side
of \eqref{12.7.4}, means that  

\begin{equation} 
                         \label{12.16.10}
 P\big(x+x_{\eta R^{2}}\in B_{\kappa\rho_{0} R },\sup_{s\leq \eta R^{2}}|x+x_{s}|<   R  \big)\geq \mu,
\end{equation}
whenever $R\in(0,\infty)$ and  
$x\in B_{\kappa R }$. For $n=1,2,...$ 
introduce ($t_{0}:=0$)
$$
R_{n}=\rho_{0}^{n-1}=\rho_{0}R_{n-1},\quad s_{n}= \eta R_{n}^{2}=\eta \rho_{0}^{2(n-1)},\quad
t_{n}=\sum_{k=1}^{n}s_{k},
$$
$$
A_{n}=\{\sup_{s\leq s_{n}}|x+x_{s+t_{n-1}}|< R_{n}\},\quad
\Pi_{n}=\bigcap_{k=1}^{n}A_{k} 
$$
and observe that by the conditional  version of \eqref{12.16.10}
on the set $\{y:=x+x_{t_{n-1}}\in B_{\kappa R_{n}}\}$ we have (a.s.)
\begin{equation}
                                                 \label{12.9.6}
P\Big(y+(x_{t_{n}}-x_{t_{n-1}})\in B_{
\kappa R_{n+1}},
\sup_{s\leq s_{n}}|y+(x_{t_{n-1}+s}-x_{t_{n-1}})|<   R_{n}
\mid \cF_{t_{n-1}}\Big)\geq \mu.
\end{equation}
Furthermore, obviously, for $n\geq 2$,
$$
P^{n}:= 
P(x+x_{t_{n}}\in B_{\kappa R_{n+1 }},\Pi_{n}  )
$$
$$
\geq P(x+x_{t_{n-1}}\in B_{\kappa R_{n-1}},\Pi_{n-1},
$$
$$
x+x_{t_{n-1}}+(x_{t_{n}}-x_{t_{n-1}})\in B_{\kappa R_{n+1}},
\sup_{s\leq s_{n}}|x+x_{t_{n-1}}+(x_{t_{n-1}+s}-x_{t_{n-1}})|<  R_{n} ),
$$
which in light of \eqref{12.9.6} yields $P^{n}\geq \mu P^{n-1}$
and since for $|x|<\kappa$
we have $P^{1}\geq \mu$ by \eqref{12.16.10}, it holds that
for $|x|<\kappa$ and all $n\geq 0$ 
\begin{equation}
                                                 \label{12.16.1}
P\big(x+x_{t_{n}}\in B_{\kappa R_{n+1 }},\sup_{s\leq t_{n}}|x+x_{s}|
< 1 \big)\geq \mu^{n+1}.
\end{equation}

Now it is convenient to consider $\eta$ as a variable
and $\kappa,d,\delta,\|\gb\|$ as fixed parameters
and not to include them in the arguments of some functions
which appear below.
Observe that
$$
t_{n}=\eta\sum_{k=0}^{n-1}\rho^{2k}_{0}(\eta)
$$
is a strictly increasing function of $\eta$ and $t_{n}\geq\eta$.
Therefore, for fixed $\eta'\in(0,1)$ and
each $n=1,2,...$ there is $\eta=\eta(n)=\eta(n,\eta')\in(0,1)$ such that
$$
\eta'=\eta\sum_{k=0}^{n-1}\rho^{2k}_{0}(\eta)\quad(=t_{n}).
$$
Clearly, the sequence $\eta(n)$ is decreasing and its limit
$\bar\eta$
is a function of $\eta'$, which is strictly positive.

Then take $\rho\in(0,\kappa \rho_{0} (\bar\eta))$ 
and define $n(\rho)=n(\rho,\eta')$ as the biggest $n\geq 1$ satisfying  
\begin{equation}
                                                 \label{12.16.2}
\kappa \rho_{0}^{n}(\bar\eta)\geq\rho
\end{equation}
that is
$$
n(\rho)=\Big\lfloor
\frac{\ln(\rho/\kappa)}{\ln \rho_{0}(\bar\eta)}\Big\rfloor.
$$

With so defined $n=n(\rho)$ in light of \eqref{12.16.2}
we have $\kappa\rho_{0}^{n}(\eta(n))\geq \rho$ and \eqref{12.16.1}
yields
$$
P\big(x+x_{\eta'}\in B_{\rho},\sup_{s\leq \eta'}|x+x_{s}|
< 1 \big)\geq \mu^{n(\rho)+1}
$$
if $\rho\in(0,\kappa \rho_{0} (\bar\eta))$ and $|x|<\kappa$.
Here
$$
\mu^{n(\rho)+1}\geq \mu \exp\Big(
\frac{\ln(\rho/\kappa)}{\ln \rho_{0}(\bar\eta)}\ln\mu\Big)
=N\rho^{\nu},
$$
where $N$ and $\nu$ are defined by the above equality.
It follows that \eqref{12.14.2}
holds with $R=1$ if $\rho\in(0,\kappa \rho_{0} (\bar\eta))$.
Then it automatically holds for larger $\rho$ 
perhaps with a different $N$. Arbitrary $R$ are
treated by self-similarity.
The lemma is proved.   \qed

{\bf Proof of Theorem \ref{theorem 12.7.2}}.
Let $R_{1}=(1-\kappa)R$ and note that 
 $\xi:= t/R_{1}^{2}-\eta$ satisfies
$$
\eta^{-1}(1-\kappa)^{-2}>\xi\geq \eta\big[(1-\kappa)^{-2}-1\big].
$$

By the conditional  version of
Lemma \ref{lemma 12.15.10} on the set
$\{z:= x+x_{\xi R_{1}^{2}}\in B_{\kappa R_{1}}(y)\}$ we have (a.s.)
$$
NP\Big(\sup_{s\in [\xi R_{1}^{2},\xi R_{1}^{2}+\eta R_{1}^{2}]}
|z+x_{s}-x_{\xi R_{1}^{2}}-y|<R_{1},
$$
$$
x+x_{\xi R_{1}^{2}+\eta R^{2}_{1}}\in B_{\rho R_{1}}(y)
\mid \cF_{\xi R_{1}^{2}}\Big)
\geq \rho^{\nu}.
$$

By Lemma \ref{lemma 12.15.1}, where we take $\rho_{0}=R_{1}/R$
and replace $\xi$ there with $\xi(1-\kappa)^{2}$,
$$
P(\sup_{s\leq \xi R_{1}^{2}}|x+x_{s}|< R,
x+x_{\xi R_{1}^{2}}\in B_{ \kappa R_{1} }(y) 
 )\geq \mu.
$$
 
By combining these two facts
and using   that
$\xi R^{2}_{1}+\eta R_{1}^{2}=t$, we obviously come
to \eqref{12.9.1}. The theorem is proved. \qed

 {\bf Acknowledgment}. The author thanks T. Yastrzhembskiy
for pointing out several mistakes and misprints in the
first draft of the paper.

\end{document}